\documentclass[final]{siamltex}         

\usepackage{graphicx,bm,amssymb,amsmath} 
\usepackage{afterpage}

\newcommand{\bi}{\begin{itemize}}
\newcommand{\ei}{\end{itemize}}
\newcommand{\ben}{\begin{enumerate}}
\newcommand{\een}{\end{enumerate}}
\newcommand{\be}{\begin{equation}}
\newcommand{\ee}{\end{equation}}
\newcommand{\bea}{\begin{eqnarray}} 
\newcommand{\eea}{\end{eqnarray}}
\newcommand{\ba}{\begin{align}} 
\newcommand{\ea}{\end{align}}
\newcommand{\bse}{\begin{subequations}} 
\newcommand{\ese}{\end{subequations}}
\newcommand{\bc}{\begin{center}}
\newcommand{\ec}{\end{center}}
\newcommand{\bfi}{\begin{figure}}
\newcommand{\efi}{\end{figure}}
\newcommand{\ca}[2]{\caption{#1 \label{#2}}}
\newcommand{\ig}[2]{\includegraphics[#1]{#2}}
\newcommand{\bmp}[1]{\begin{minipage}{#1}}
\newcommand{\emp}{\end{minipage}}

\newcommand{\tbox}[1]{{\mbox{\tiny #1}}}
\newcommand{\mbf}[1]{{\mathbf #1}}
\newcommand{\half}{\mbox{\small $\frac{1}{2}$}}
\newcommand{\RR}{\mathbb{R}}
\newcommand{\ZZ}{\mathbb{Z}}
\newcommand{\eps}{\varepsilon}
\newcommand{\bigO}{{\mathcal O}}
\newcommand{\intR}{\int_{-\infty}^\infty}
\newcommand{\eg}{e.g.\ }
\newcommand{\ie}{i.e.\ }
\newcommand{\etal}{et al.\ }

\newtheorem{thm}{Theorem}

\newtheorem{rmk}[thm]{Remark}
\newcommand{\xx}{\mbf{x}}
\newcommand{\sss}{\mbf{s}}

\newcommand{\kk}{\mbf{k}}
\newcommand{\KK}{{\mathcal I}}     
\newcommand{\freq}{\beta}          
\newcommand{\rat}{\sigma}          
\newcommand{\ppsi}{{\tilde\psi}}   
\newcommand{\rmax}{r_\tbox{dyn}}    
\newcommand{\al}{\alpha}           
\newcommand{\NU}{{nonuniform}}       
\newcommand{\U}{{uniform}}
\newcommand{\KB}{Kaiser--Bessel}
\newcommand{\FT}{Fourier transform}

\begin{document}

\title{A parallel non-uniform fast Fourier transform library
  based on an ``exponential of semicircle'' kernel}

\author{Alex Barnett%
  \thanks{Flatiron Institute, Simons Foundation, New York, NY, USA}
  \and
  Jeremy Magland%
  \thanks{Flatiron Institute, Simons Foundation, New York, NY, USA}
  \and
  Ludvig af Klinteberg%
  \thanks{Department of Mathematics, Simon Fraser University, Burnaby, BC, Canada}
  }
  \date{\today}
\maketitle
\begin{abstract}
  The nonuniform fast Fourier transform
  (NUFFT) generalizes the FFT to off-grid data.
  %
  Its many applications
  include image reconstruction, data analysis,
  and the numerical solution of differential equations.
  %
  We present FINUFFT, an efficient parallel library
  for type 1 (\NU\ to \U), type 2 (\U\ to \NU), or
  type 3 (\NU\ to \NU) transforms, in dimensions 1, 2, or 3.
  It uses minimal RAM, requires no precomputation or plan steps,
  and has a simple interface to several languages.
  %
  We perform the expensive 
  spreading/interpolation between nonuniform points and the fine grid
  via a simple new kernel---the
  ``exponential of semicircle''
  $e^{\beta \sqrt{1-x^2}}$ in $x\in[-1,1]$---%
  in a cache-aware load-balanced multithreaded implementation.
  The deconvolution step requires the Fourier transform of the kernel,
  for which we propose efficient numerical quadrature.
  For types 1 and 2,
  rigorous error bounds asymptotic in the kernel width
  approach the fastest known exponential rate, namely
  that of the Kaiser--Bessel kernel.
  We benchmark 
  against several popular CPU-based libraries,
  showing favorable speed and memory footprint,
  especially in three dimensions when high accuracy and/or
  clustered point distributions are desired.
\end{abstract}

\section{Introduction}


The need for fast algorithms for spectral analysis of
non-uniformly sampled data arose
soon after the popularization of the FFT in the 1960s.
Many early methods came from
signal processing \cite{oppenheim71} and
astronomy \cite{thompson74,meisel78,rybicki} \cite[Sec.~13.8]{numrec},
but it was not until
the 1990s that Dutt--Rokhlin \cite{dutt}
gave the first rigorous analysis of a convergent scheme.
The NUFFT has since become crucial in many areas of science and engineering.
Several imaging methods, including magnetic resonance imaging
(MRI) \cite{suttonfield,nufft3}, X-ray computed tomography (CT) \cite{fourmont},
ultrasound diffraction tomography \cite{ultranufft},
and synthetic aperture radar \cite{andersson12},  
sample the Fourier transform at non-Cartesian points \cite{gelbrecon,fastsinc}
hence require the NUFFT or its inverse for accurate reconstruction.
Real-time Fourier-domain optical coherence tomography (OCT)
relies on rapid one-dimensional NUFFTs \cite{octnufft}.
Periodic electrostatic and Stokes problems
are commonly solved by fast
``particle-mesh Ewald'' summation, whose spectral part is equivalent to
a pair of NUFFTs \cite{Lindbo11,nestlerPME}.
Spectrally-accurate function interpolation may be efficiently performed with the
NUFFT \cite[Sec.~6]{usingnfft} \cite{gimbutasgrid}.
The numerical approximation of Fourier transforms using
non-Cartesian or adaptive
quadrature grids arises in heat solvers \cite{nufft3},
cryo electron microscopy \cite{steerablePCA,cryo},
and electromagnetics \cite{emnufft}.
Many more applications are found in the reviews
\cite{warenufft,kunisthesis,nfftchap,usingnfft,nufft}.

Our purpose is to describe and benchmark a 
general-purpose software library
for the NUFFT that achieves high efficiency with an open-source compiler
by combining mathematical and implementation innovations.
The computational task is
to approximate, to a requested relative accuracy $\eps$,
the following exponential sums.
Let $d=1,2$ or 3 be the spatial dimension.
Let $N_i$ be the number of desired Fourier modes in dimension $i = 1,\dots,d$;
in each dimension the Fourier mode (frequency) indices are
$$
\KK_{N_i} := \left\{
\begin{array}{ll} \{-N_i/2,\ldots,N_i/2-1\}, & N_i \mbox{ even},\\
\{-(N_i-1)/2,\ldots,(N_i-1)/2\}, & N_i \mbox{ odd}~.
\end{array}\right.
$$
The full set of mode indices is the Cartesian product that we denote by
$$
\KK = \KK_{N_1,\dots,N_d} := \KK_{N_1} \times \dots \times \KK_{N_d}~,
$$
containing a total number of modes $N=N_1\cdots N_d$.
The $M$ \NU\ points have locations $\xx_j$, $j=1,\ldots,M$,
which may be taken to lie in $[-\pi,\pi)^d$,
with corresponding strengths $c_j\in\mathbb{C}$.
Then the type 1 NUFFT
(also known as the ``adjoint NFFT'' \cite{nfftchap,usingnfft})
approximates the outputs%
\footnote{Note that our normalization differs from that of \cite{dutt,nufft}.}
\be
f_\kk :=
\sum_{j=1}^M c_j e^{i \kk\cdot \xx_j}~,
\quad 
\kk \in \KK~.
\qquad \mbox{(Type 1, \NU\ to \U)}
\label{1}
\ee
This may be interpreted as computing, for the $2\pi$-periodic box,
the $N$ Fourier series coefficients of the distribution
\be
f(\xx) := \sum_{j=1}^M c_j \delta(\xx-\xx_j) ~.
\label{fdist}
\ee
Up to normalization, \eqref{1} generalizes the discrete Fourier transform (DFT),
which is simply the uniform case, \eg in 1D, $\xx_j = 2\pi j/M$ with $M=N_1$.

The type 2 transform (or ``NFFT'') is the adjoint of the type 1.
Unlike in the DFT case, it is not generally
related to the inverse of the type 1.
It evaluates the Fourier series with given coefficients
$f_\kk$ at arbitrary 
target points $\xx_j$.
  That is,
  \be
  c_j := \sum_{\kk\in\KK} f_\kk e^{-i \kk\cdot \xx_j}~,
  \quad j=1,\dots, M~.
\qquad \mbox{(Type 2, \U\ to \NU)}
\label{2}
\ee
Finally, the more general type 3 transform
\cite{nufft3} (or ``NNFFT'' \cite{usingnfft})
may be interpreted as evaluating the
Fourier {\em transform} of the non-periodic distribution \eqref{fdist}
at arbitrary target frequencies
$\sss_k$ in $\RR^d$, $k=1,\dots, N$, where $k$ is a plain integer index.
That is,
\be
f_k := \sum_{j=1}^M c_j e^{i \sss_k \cdot \xx_j}~,
  \quad k=1,\dots, N~.
\qquad \mbox{(Type 3, \NU\ to \NU)}
\label{3}
\ee
All three types of transform \eqref{1}, \eqref{2} and \eqref{3}
consist simply of computing exponential sums
that 
naively require $\bigO(NM)$ work.
NUFFT algorithms compute these
sums, to a user-specified relative tolerance $\eps$,
in close to linear time in $N$ and $M$.

\begin{rmk}
In certain settings the above sums may be interpreted as quadrature formulae
applied to evaluating a Fourier transform of a function.
However, these tasks are not to be confused with the ``inverse NUFFT''
(see problems 4 and 5 in \cite{dutt,nufft})
which involves,
for instance, treating \eqref{2} as a linear system to be solved for
$\{f_\kk\}$, given the right-hand side $\{c_j\}$.
For some \NU\ distributions this linear system can be
very ill-conditioned.
This inverse NUFFT is common in Fourier imaging
applications; a popular solution method is to use
conjugate gradients to solve
the preconditioned normal equations, exploting
repeated NUFFTs for the needed matrix-vector multiplications
\cite{fessler,fourmont,fastsinc,gelbrecon} \cite[Sec.~3.3]{townsendnufft}.
Thus, efficiency gains reported here would also accelerate this
inversion method.
See \cite{duttcotfmm,kircheis} for other approaches.
We will not explicitly address the inversion problem here.
\end{rmk}

\bfi[t]  
\ig{width=\textwidth}{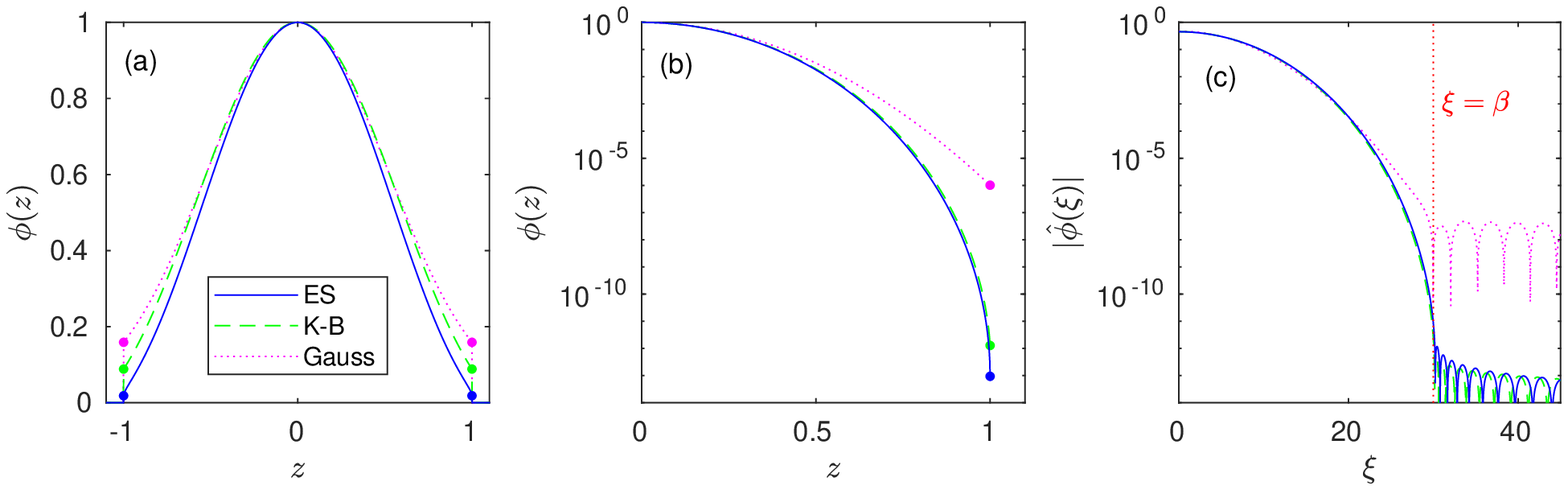}
\ca{The proposed ES spreading kernel
  \eqref{ES} (solid blue lines); the \KB\ kernel \eqref{KB}
  (dashed green); and the best truncated Gaussian (dotted pink)
  $\phi(z) = e^{-0.46 \freq z^2}$ in $|z|\le 1$. 
  (a) shows all kernels for $\freq=4$. The discontinuities at
  $\pm 1$ are highlighted by dots.
  (b) shows a logarithmic plot (for positive $z$) of the
  kernels for $\freq=30$
  (corresponding to a spreading width of $w=13$ grid points).
  The graph for ES is a quarter-circle.
  (c) shows the magnitude of the kernel Fourier transforms,
  and the ``cutoff'' frequency $\xi=\freq$.
  ES and KB have shape close to a quarter-ellipse in $|\xi|<\freq$
  (see \eqref{EShat1} and \eqref{KBhat}).
  All have exponentially small values for $|\xi|>\freq$, but
  the Gaussian has exponential convergence rate in terms of $\freq$
  only around half that of ES or KB.
}{f:kernel}
\efi

\subsection{Prior algorithms, kernels, and implementations}
\label{s:prior}

There are two main approaches to the fast approximation
of the sums \eqref{1} or \eqref{2}, both of which build upon the FFT:
1) Interpolation between \NU\ points and an upsampled regular grid,
combined with an upsampled FFT and correction in Fourier space
\cite{dutt,beylkinnufft,fessler,nufft,usingnfft};
or
2) Interpolation to/from an $N$-point (\ie not upsampled)
regular grid, combined with the $N$-point FFT.
In the univariate (1D) case, there are several variants of the second approach:
Dutt--Rokhlin \cite{duttcotfmm} proposed
spectral Lagrange interpolation (using the cotangent kernel
applied via the fast multipole method), combined with a
single FFT.
Recently, Ruiz-Antol\'in and Townsend \cite{townsendnufft}
proposed a stable Chebyshev approximation in intervals centered about
each \U\ point,
which needs an independent $N$-point FFT for each of the
$\bigO(\log 1/\eps)$ coefficients, but is embarrassingly
parallelizable. This improves upon earlier work \cite{anderson96}
using Taylor approximation that was numerically unstable
without upsampling \cite[Ex. 3.10]{kunisthesis} \cite{townsendnufft}.

We now turn to the first, and most popular,
of the two above approaches.
For the type 1 and type 2 transforms
one sets up a regular fine grid of $n = \sigma^d N$ points
where the upsampling factor in each dimension, $\sigma>1$, is a small constant
(typically $\le 2$).
Taking the type 1 as an example, there are three steps.
Step 1 evaluates on the fine grid the convolution of \eqref{fdist}
with a smooth kernel 
function $\psi$,
whose support has width $w$ fine grid points in each dimension
(see Fig.~\ref{f:spreadalias}(a)).
This ``spreading'' requires $w^d M$ kernel evaluations.
Step 2 applies the FFT on the $n$-point grid, needing $\bigO(N \log N)$ work.
Step 3 extracts the lowest $N$ frequencies from the result,
then divides by $\hat \psi$, the \FT\ of the kernel,
evaluated at each of these frequencies;
this is called deconvolution or roll-off correction.
There is a class of kernels, including all those we discuss below,
whose analysis gives an error $\eps$ decreasing
exponentially (up to weak algebraic prefactors) with $w$,
hence one may choose $w \approx c |\log\eps|$.
Thus the total effort for the NUFFT is $\bigO(M |\log \eps|^d + N \log N)$.

The choice of spreading kernel $\psi$ has a fascinating history.
A variety of kernels were originally used for
``gridding'' in the imaging community
(\eg see \cite{jackson91,nfftchap,usingnfft}).
The truncated Gaussian kernel (see Fig.~\ref{f:kernel})
was the first for which an exponential
convergence rate with respect to $w$ was shown \cite{dutt}.
This rate was improved by Steidl \cite{steidl98,elbel}:
fixing $\sigma$, for an optimally chosen Gaussian width parameter
the error is $\eps = \bigO(e^{-\half\pi(1-(2\sigma-1)^{-1}) w})$.
Beylkin \cite{beylkinnufft} proposed B-splines for $\psi$, with the
estimate $\eps = \bigO((2\sigma-1)^{-w})$.
In both cases, it is clear that
increasing $\sigma$ improves the convergence rate;
however, since the cost of the upsampled FFT grows at least like $\sigma^d$,
a tradeoff arises.
In practice, many studies have settled on $\sigma=2$ for general use
\cite{jackson91,fourmont,fessler,usingnfft,nufft,ou}.
For this choice, both the above Gaussian
and B-spline rates imply that
$|\log_{10} \eps|$, the number of correct digits, is
approximately $0.5 w$.
For instance, to achieve 12 digits, a spreading width $w=24$ is needed
\cite[Remark~2]{nufft}.

However, Jackson et al.\ \cite{jackson91}
realized that the criteria for a good kernel $\psi$ are very similar
to those for a good {\em window function} in digital signal processing (DSP).
To summarize these criteria,
\bi
\item[(a)] The numerical support of $\psi$ in fine grid points, $w$,
  should be as small as possible, in order
  to reduce the $\bigO(w^dM)$ spreading cost.
\item[(b)] A certain norm of $\hat\psi(k)$ in the ``tails''
  $|k|\ge (\sigma-\half)N$
  should be 
  as small as possible, relative to values in the central range $|k|<N/2$;
  see Fig.~\ref{f:spreadalias}(b).
  \ei
The two criteria conflict:
(a) states that $\psi$ should be narrow, but (b), which
derives from {\em aliasing error}, implies that $\psi$ should
be smooth.
(We postpone the rigorous statement of (b) until \eqref{epsest}.)
Is has been known since the work of Slepian and coworkers in
the 1960s \cite{slepian65}
that, if one chose $L^2$-norms
in (b), the family of prolate spheroidal wavefunctions (PSWF) of order zero
\cite{osipov}
would optimize the above criteria.
It was also DSP folklore  \cite{kaiser} that the 
``\KB'' (KB) kernel,
\be
\phi_{\tbox{KB},\freq}(z) := \left\{
\begin{array}{ll}I_0(\freq\sqrt{1-z^2}) / I_0(\freq)~, & |z|\le 1\\
  0~,& \mbox{otherwise}\end{array}\right.
\label{KB}
\ee
scaled here to have support $[-1,1]$,
where $I_0$ is the regular modified Bessel function of order zero
\cite[(10.25.2)]{dlmf},
well approximates the PSWF.
However, unlike the PSWF, which is tricky to evaluate accurately
\cite{osipov}, \eqref{KB} needs only
standard special function libraries \cite{CEPHES}.
Its \FT\ (using the convention \eqref{FT}) is known analytically%
\footnote{This pair appears to be a discovery of B. F. Logan,
and its use pioneered in DSP by J. F. Kaiser, both at
Bell Labs, in the 1960s \cite{kaiserinterview}.
Curiously, the pair seems absent from all standard tables \cite[\S 6.677]{GS8}
\cite[\S 2.5.25]{prudnikov1} \cite[\S 2.5.10]{prudnikov2}.}
\cite{kaiser},
\be
\hat\phi_{\tbox{KB},\freq}(\xi) = \frac{2}{I_0(\freq)}
\frac{\sinh \sqrt{\freq^2-\xi^2}}{\sqrt{\freq^2-\xi^2}}
~.
\label{KBhat}
\ee
This transform pair \eqref{KB}--\eqref{KBhat} is plotted in green in
Fig.~\ref{f:kernel}.

Starting with imaging applications in the 1990s,
the KB kernel \eqref{KB} was introduced for the NUFFT
\cite{jackson91,nfftchap,fessler}.
Note that the function \eqref{KBhat}, truncated to $[-\freq,\freq]$,
outside of which it is exponentially small,
may instead be used as the spreading kernel \cite{fourmont,usingnfft}.
This latter approach---which we call ``backward KB''%
\footnote{The distinction between forward and backward use of
  the KB pair is unclear in the literature.}%
---%
has the computational advantage of spreading with
cheaper sinh rather than $I_0$ evaluations.
The error analyses of the two variants turn out to be equivalent.
Despite being only conditionally convergent,
the tail sum of \eqref{KBhat} needed for criterion (b)
may be bounded rigorously;
this subtle analysis is due to Fourmont \cite[p.~30--38]{fourmontthesis}
\cite[Sec.~4]{fourmont}.
Its optimal convergence with $w$, summarized in \cite[p.~30-31]{pottshabil}
\cite[App.~C]{nfft}, is (see \eqref{epsest} for the definition of
the error $\eps_\infty$),
\be
\eps_\infty \;\le\; 4\pi   (1-1/\rat)^{1/4}
\left(\sqrt{\frac{w-1}{2}}+\frac{w-1}{2}\right)
e^{-\pi(w-1)\sqrt{1-1/\rat}}
~.
\label{KBerr}
\ee
This is the fastest known exponential error rate of any kernel,
equalling that of the PSWF \cite{nufftanal}:
for the choice $\sigma=2$ gives
over $0.9w$ correct digits.
This is {\em nearly twice} that of the Gaussian;
12 digits are reached with only $w=13$.
Attempts to further optimize this kernel
give only marginal gains \cite{fessler},
unless restricted to cases with specific decay of
the mode data $f_\kk$ \cite{nestler},
or minimal upsampling ($\rat\approx 1$) \cite{l2jacob}.

Turning to software implementations, most
are based upon the Gaussian or KB kernels (in both its variants).
Greengard--Lee \cite{nufft} presented ``fast Gaussian gridding''
which reduced the number of exponential function evaluations
from $w^dM$ to $(d+1)M$, resulting in a several-fold acceleration of
the spreading step.
This was implemented by those authors in
a general-purpose single-threaded CMCL Fortran library \cite{cmcl}.
The mature general-purpose
NFFT code of Keiner--Kunis--Potts \cite{nfft,usingnfft}
is multithreaded \cite{volkmer}
and uses backward KB by default (although fast Gaussian gridding is available).
It allows various {\em precomputations}
of kernel values (requiring a ``plan'' stage),
demanding a larger RAM footprint
but accelerating repeated calls with the same points.
There are also several codes specialized to MRI,
including MIRT (which uses full precomputation of the KB kernel)
by Fessler \etal \cite{MIRT}, and recently BART \cite{BART}
and pynufft \cite{pynufft}.
Various specialized GPU implementations also exist (reviewed in \cite{ou}),
mostly for MRI \cite{cunfft,gpunufft} or OCT \cite{octnufft}.
Unlike general-purpose codes, these specialized packages tend
to have limited accuracy or dimensionality, and tend not to document
precisely what they compute.

\subsection{Contribution of this work}

We present a general purpose documented CPU-based multithreaded C++ library
(FINUFFT) \cite{finufftgit}
that is efficient without needing any precomputation stage.
This means that the RAM overhead is
very small
and the interface simple.
For medium and large problems in 2D and 3D its speed is
competitive with state-of-the-art CPU-based codes.
In some cases, at high accuracies,
FINUFFT is faster than all known CPU-based codes by a factor of 10.
The packages against which we benchmark are
listed in Table~\ref{t:codes}.

We spread with a new
``exponential of semicircle'' (ES)
kernel (see Fig.~\ref{f:kernel}),
\be
\phi_\freq(z) :=
\left\{\begin{array}{ll}
e^{\freq (\sqrt{1-z^2}-1)}, & |z|\le 1\\
0, & \mbox{otherwise}
\end{array}
\right.
\label{ES}
\ee
which has error convergence rate
arbitrarily close
to that of \eqref{KBerr}; see Theorem~\ref{t:ESerr}.
It is simpler and faster to evaluate than either of the KB pair
\eqref{KB}--\eqref{KBhat}, yet has essentially identical error.
We demonstrate further acceleration via piecewise polynomial
approximation.
\eqref{ES} has no known analytic \FT, yet we
can use {\em numerical quadrature} to evaluate $\hat\phi_\freq$ when needed,
with negligible extra cost.
Unlike interpolation from the fine grid (needed in type 2),
spreading (needed for type 1)
does not naturally parallelize over \NU\ points, because
of collisions between threads writing to the output grid.
However, we achieve efficiency in this case
by adaptively blocking into auxiliary fine grids, after bin-sorting the points.

The rest of the paper is structured as follows.
The next section outlines the software interfaces.
In section~\ref{s:alg} we describe the
algorithms and parameter choices in full,
including various novelties in terms
of quadrature and type 3 optimization.
In Section~\ref{s:err} we
summarize a rigorous aliasing error bound for the ES kernel,
and use this to justify the choice of $w$ and $\freq$.
We also explain the gap between this bound and empirically observed
relative errors.
Section~\ref{s:perf}
compares the speed and accuracy performance against other libraries,
in dimensions 1, 2, and 3.
We conclude in section~\ref{s:conc}.

\section{Use of the FINUFFT library}



The basic interfaces are very simple \cite{finufftgit}.
From C++, with {\tt x} a {\tt double} array of {\tt M} source points,
{\tt c} a complex ({\tt std::complex<double>}) array of {\tt M} strengths,
and {\tt N} an integer number of desired output modes,
\vspace{.5ex}
\begin{verbatim}
  status = finufft1d1(M,x,c,isign,tol,N,f,opts);
\end{verbatim}
\vspace{.5ex}
computes the 1D type 1 NUFFT with relative precision {\tt tol}
(see section~\ref{s:params}), writing the
answer into the preallocated complex array {\tt f},
and returning zero if successful.
Setting {\tt isign} either $1$ or $-1$ controls the
sign of the imaginary unit in \eqref{1}.
{\tt opts} is a struct
defined by the required header {\tt finufft.h}
and initialized by \verb+finufft_default_options(&opts)+,
controlling various options.
For example, setting {\tt opts.debug=1} prints internal timings,
whereas {\tt opts.chkbnds=1} includes an initial check whether all
points lie in the valid input range $[-3\pi,3\pi]$.
The above is one of nine routines with similar interfaces (types 1, 2, and 3,
in dimensions 1, 2, and 3).
The code is lightweight, at around 3300 lines of C++ (excluding
interfaces to other languages).
DFTs are performed by FFTW \cite{fftw}, which is the only dependency.

Interfaces from C and Fortran are similar to the above,
and require linking with \verb?-lstdc++?.
From high-level languages one may call,
\vspace{.5ex}
\begin{verbatim}
  [f status] = finufft1d1(x,c,isign,tol,N,opts);    % MATLAB or octave
  status = finufftpy.nufft1d1(x,c,isign,tol,N,f)    # python (numpy)
\end{verbatim}
\vspace{.5ex}
Here $M$ is inferred from the input sizes.
There also exists a julia interface \cite{finufft.jl}.

\begin{rmk}
  The above interface, since it does not involve any ``plan'' stage,
  incurs a penalty for repeated small problems ($N$ and $M$
  of order $10^4$ or less), traceable
  to the overhead (around 100 microseconds per thread in our tests)
  for calling \verb+fftw_plan()+ present when FFTW reuses stored wisdom.
  To provide maximal throughput for repeated small problems (which are yet not
  small enough that a ZGEMM matrix-matrix multiplication approach wins),
  we are adding interfaces that handle multiple inputs or allow a plan stage.
  At the time of writing these are available in 2D only, and will be
  extended in a future release.
  \label{r:small}
\end{rmk}

\section{Algorithms}
\label{s:alg}

For type 1 we use the standard three-step procedure sketched above
in Section~\ref{s:prior}. For type 2 the steps are reversed.
Type 3 involves a combination of types 1 and 2.
Our Fourier transform convention is
\be
\hat\phi(k) = \intR \phi(x) e^{ikx} dx
~,\qquad
\phi(x) = \frac{1}{2\pi} \intR \hat\phi(k) e^{-ikx} dx
~.
\label{FT}
\ee
For the default upsampling factor $\rat=2$,
given the requested relative tolerance $\eps$,
the kernel width $w$ and ES parameter $\freq$ in \eqref{ES} are set via
\be
w \;=\; \lceil \log_{10} 1/\eps \rceil + 1 ~,
\qquad
\freq \;=\; 2.30 \, w ~.
\label{params}
\ee
The first formula may be
summarized as:
the kernel width is one more than the desired number of accurate digits.
We justify these choices in section~\ref{s:params}.
(FINUFFT also provides a low-upsampling option $\rat=5/4$,
which is not tested in this paper.)


\bfi[t]  
\ig{width=5.2in}{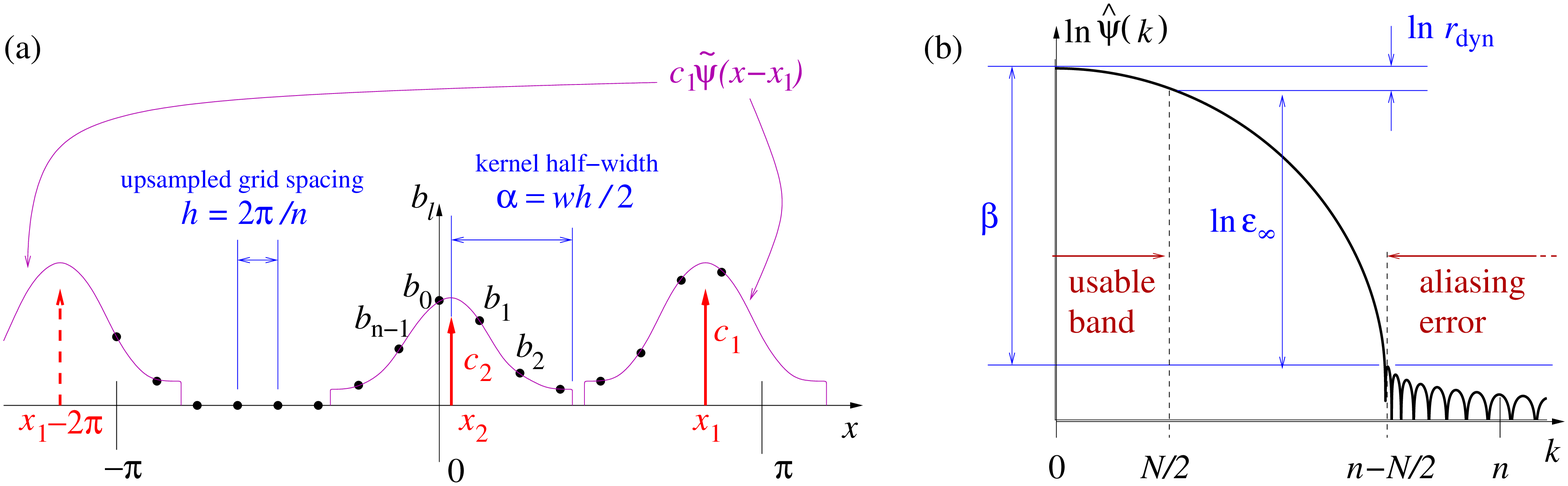}
\ca{(a) 1D illustration of spreading from \NU\ points to the grid
  values $b_l$, $l=0,\dots,n-1$ (shown as dots) needed for type 1.
  For clarity, only two \NU\ points $x_1$ and $x_2$ are shown;
  the former results in periodic wrapping of the effect of the kernel.
  (b) Semi-logarithmic plot of the (positive half of the)
  \FT\ of the rescaled
  kernel $\psi(x)$, showing the usable frequency domain (and the
  dynamic range $\rmax$ over this domain), and
  a useful approximate relationship between the
  aliasing error bound $\eps_\infty$
  and the ES kernel 
  parameter $\beta$.
  From \eqref{EShat1}, below cutoff
  the curve is well approximated by
  a quarter 
  ellipse.
}{f:spreadalias}
\efi

\subsection{Type 1: \NU\ to \U}
\label{s:type1}

We describe the algorithm to compute $\tilde f_\kk$, an approximation
to the exact $f_\kk$ defined by \eqref{1}.

\subsubsection{1D case}
\label{s:1d1}

We use $x_j$ to denote \NU\ source points, and $k\in\KK$
to label the $N=N_1$ output modes.
For FFT efficiency the
DFT size $n$ is chosen to be the smallest integer of the form
$2^q3^p5^r$ not less than $\rat N$ nor $2w$, the latter condition
simplifying the spreading code.

{\bf Step 1 (spreading).}
From now on we abbreviate the ES kernel $\phi_\freq$ in \eqref{ES} by $\phi$.
We rescale the kernel so that its support becomes
$[-\al,\al]$, with
\be
\alpha := wh/2 = \pi w/n
~,
\label{al}
\ee
where $h := 2\pi/n$ is the upsampled grid spacing.
This rescaled kernel is denoted
\be
\psi(x) := \phi(x/\al)~,
\qquad \mbox{thus} \quad \hat\psi(k) = \al \hat\phi(\al k)~,
\qquad \mbox{(1D case)}
\label{psi1}
\ee
and its periodization is
\be
\ppsi(x) := \sum_{m\in\ZZ} \psi(x-2\pi m)
~.
\qquad \mbox{(1D case)}
\label{ppsi1}
\ee
We then compute, at a cost of $wM$ kernel evaluations, 
the periodic discrete convolution
\be
b_l = \sum_{j=1}^M c_j \ppsi(lh - x_j)
~, \qquad \mbox{for } l=0,\dots,n-1
~,
\label{bl1}
\ee
as sketched in Fig.~\ref{f:spreadalias}(a).

{\bf Step 2 (FFT).}
We use the FFT to evalute the $n$-point DFT
\be
\hat{b}_k = \sum_{l=0}^{n-1} e^{2\pi i lk/n} b_l ~, \qquad \mbox{ for } k\in\KK
~.
\label{dft1}
\ee
Note that the output index set $\KK$ is cyclically equivalent to the
usual FFT index set $k=0,\dots,n-1$.

{\bf Step 3 (correction).}
We truncate to the central $N$ frequencies, then
diagonally scale (deconvolve) the amplitudes array, to give the outputs
\be
\tilde f_k = p_k \hat{b}_k ~, \qquad \mbox{ for } k\in\KK
~,
\label{pb1}
\ee
where a good choice of the correction factors $p_k$ comes from
samples of the kernel Fourier transform,%
\footnote{It is tempting instead to set $p_k$ to be the DFT 
  of the grid samples of the kernel $\{\ppsi(lh)\}_{l=0}^{n-1}$.
  However, in our
  experience this causes around twice the error of \eqref{pk1},
  as can be justified by the discussion in section~\ref{s:err}.
Fessler and Sutton \cite[Sec.~V.C.3]{fessler} report a similar finding.}
\be p_k = h / \hat\psi(k)~ = 2/(w\hat\phi(\alpha k)), \qquad k\in\KK
\label{pk1}
~.
\ee
A new feature of our approach is that
we approximate $\hat\psi(k)$
by applying Gauss--Legendre quadrature
to the Fourier integral, as follows.
This allows kernels without an analytically known Fourier transform
to be used without loss of efficiency.
Let $q_j$ and $w_j$ be the nodes and weights
for a $2p$-node quadrature on $[-1,1]$.
Since $\phi$ is real and even, only the $p$ positive nodes are needed, thus,
\be
\hat\psi(k) \;=\; 
\int_{-\al}^{\al} \psi(x) e^{ikx} dx
\;\approx\;
wh \sum_{j=1}^p w_j \phi(q_j) \cos (\al k q_j)
~.
\label{quad}
\ee
By a convergence study,
we find that $p\ge 1.5 w+2$ (thus a maximum quadrature spacing close to $h$)
gives errors less than $\eps$, over the needed range $|k|\le N/2$.
A rigorous quadrature bound would be difficult due to the small square-root
singularities at the endpoints in \eqref{ES}.
The cost of the evaluation of $p_k$ is $\bigO(pN)$,
and naively would involve $pN$ cosines.
By exploiting the fact that, for each quadrature point $q_j$,
successive values of $e^{i \al k q_j}$ over the regular $k$ grid are
related by a constant phase factor, these cosines
can be replaced by only $p$ complex exponentials, and $pN$ adds and multiplies,
giving an order of magnitude acceleration.
We call this standard 
trick
``phase winding.''\footnote{In the code, see the function
  {\tt onedim\_fseries\_kernel} in {\tt src/common.cpp}}

\subsubsection{The case of higher dimension $d>1$}

In general, different
fine grid sizes are needed in each dimension.
We use the same recipe, so that $n_i \ge \rat N_i$, $n_i \ge 2w$,
$n_i = 2^{q_i}3^{p_i}5^{r_i}$, $i=1,\dots,d$.
The kernel is a periodized product of scaled 1D kernels,
\be
\psi(\xx) = \phi(x_1/\al_1) \cdots \phi(x_d/\al_d)~,
\qquad
\ppsi(\xx) := \sum_{\mbf{m} \in \ZZ^d} \psi(\xx - 2\pi\mbf{m})
~,
\label{ppsi}
\ee
where $\al_i=\pi w/n_i$.
Writing $h_i:=2\pi/n_i$ for the fine grid spacing in each dimension,
and $\mbf{l}:=(l_1,\dots,l_d)$ to index each fine grid point,
the discrete convolution 
becomes
\be
b_\mbf{l} = \sum_{j=1}^M c_j \ppsi((l_1h_1,\dots,l_dh_d)-\xx_j)~,
\qquad l_i=0,\dots,n_i-1, \quad i=1,\dots,d
~.
\label{bl}
\ee
In evaluating \eqref{bl}, separability
means that only $wd$ kernel evaluations are needed per source point:
the $w^d$ square or cube of $\ppsi$ values is then filled by an outer product.
The DFT \eqref{dft1} generalizes in the standard way to
multiple dimensions.
Finally, the correction factor is also separable,
\be
p_\kk = h_1\dots h_d \hat\psi(\kk)^{-1} = (2/w)^{d}
(\hat\phi(\al_1 k_1) \cdots \hat\phi(\al_d k_d))^{-1}
~, \qquad \kk \in \KK~,
\label{pk}
\ee
so that only $d$ repetitions of \eqref{quad} are needed,
followed by an outer product.

\subsection{Type 2: \U\ to \NU}
\label{s:2}

To compute $\tilde c_j$, an approximation to $c_j$ in \eqref{2},
we reverse the steps for the type 1.
Given the number of modes $N$, and the precision $\eps$,
the choices of $n$, $w$ and $\beta$ are as in the type 1.
From now on we stick to the case of general dimension $d$.

{\bf Step 1 (correction).}
The input coefficients $f_\kk$ are pre-corrected (amplified) and zero-padded
out to the size of the fine grid,
\be
\hat b_\kk = \left\{\begin{array}{ll}p_\kk f_\kk~, & \kk \in \KK \\
0~, & \kk \in \KK_{n_1,\dots,n_d} \backslash \KK
\end{array}\right.
\ee
with the same amplification factors $p_\kk$ as in \eqref{pk}.

{\bf Step 2 (FFT).}
This is just as in type 1. Writing the general dimension
case of \eqref{dft1}, with the index vectors $\mbf{l}$
and $\kk$ (and their ranges) swapped,
\be
b_\mbf{l} = \sum_{\kk\in\KK_{n_1,\dots,n_d}}
e^{2\pi i (l_1k_1/n_1 + \dots + l_dk_d/n_d)}
\,\hat b_\kk ~, \qquad \mbox{ for }
\quad l_i=0,\dots,n_i-1, \quad i=1,\dots,d
~.
\label{dft}
\ee

{\bf Step 3 (interpolation).}
The adjoint of spreading is interpolation, which
outputs a weighted admixture of the
grid values near to each target point.
The output is then
\be
\tilde c_j = \sum_{l_1=0}^{n_1-1} \cdots \sum_{l_d=0}^{n_d-1}
b_\mbf{l} \ppsi((l_1h_1,\dots,l_dh_d) - \xx_j)
~.
\label{interp}
\ee
As with the type 1,
because of separability,
this requires $wd$ evaluations of the kernel function, and $w^d$ flops,
per target point.

\subsection{Type 3: \NU\ to \NU}

This algorithm is more involved, but
is a variant of standard ones
\cite[Alg.~3]{nufft} \cite[Alg.~2]{elbel} \cite{nufft3}
\cite[Sec.~1.3]{nfftchap}.
Loosely speaking, it is
``a type 1 wrapped around a type 2,'' where
the type 2 replaces the middle FFT step of the type 1.
Given $\eps$, we choose $w$, $\freq$, and $p$ as before.
We will present the choice of $n_i$ shortly (see ``step 0'' below).
It will involve the following bounds on
source and target coordinates $\xx_j = (x_j^{(i)},\dots,x_j^{(d)})$
and $\sss_k = (s_k^{(i)},\dots,s_k^{(d)})$:
\be
X_i := \max_{j=1,\dots,M} |x_j^{(i)}|
~,\quad
S_i := \max_{k=1,\dots,N} |s_k^{(i)}|
~,\qquad\mbox{ for } i=1,\dots,d~.
\label{XS}
\ee

{\bf Step 1 (dilation and spreading).}
For spreading onto a grid on $[-\pi,\pi)^d$,
a dilation factor $\gamma_i$ needs to be chosen
for each dimension $i=1,\ldots,d$ 
such that the rescaled sources ${x'_j}^{(i)} := x_j^{(i)}/\gamma_i$
lie in $[-\pi,\pi)$. Furthermore these rescaled coordinates must be
at least $w/2$ grid points
from the ends $\pm\pi$ in order
to avoid wrap-around of mode amplitudes in Step 2.
This gives a condition relating $n_i$ and $\gamma_i$,
\be
X_i/\gamma_i \;\le\; \pi(1 - w/n_i)
~, \qquad i=1,\dots,d~.
\label{cond1}
\ee
We may then rewrite \eqref{3} as
$f_k := \sum_{j=1}^M c_j e^{i \sss'_k \cdot \xx'_j}$, $k=1,\dots, N$,
where ${s'_k}^{(i)} = \gamma_i s_k^{(i)}$.

We spread these rescaled sources
$\xx'_j=({x'_j}^{(i)},\dots,{x'_j}^{(d)})$
onto a regular grid using the usual
periodized kernel \eqref{ppsi}, to get
\be
\hat b_\mbf{l} = \sum_{j=1}^M c_j \ppsi((l_1h_1,\dots,l_dh_d)-\xx'_j)~,
\qquad \mbf{l} \in \KK_{n_1,\dots,n_d}~.
\label{blt3}
\ee
Unlike before, here we have chosen a (cyclically equivalent) output index
grid centered at the origin, because we
shall now interpret $\mbf{l}$ as a Fourier mode index.

{\bf Step 2 (Fourier series evaluation via type 2 NUFFT).}
Treating $\hat b_\mbf{l}$ from \eqref{blt3} as Fourier series coefficients, we evaluate this series at rescaled target points using the type 2 NUFFT
(see section~\ref{s:2}), thus,
\be
b_k = \sum_{\mbf{l} \in \KK_{n_1,\dots,n_d}}
\!\! \hat b_\mbf{l} \, e^{i\mbf{l}\cdot \sss''_k}
\,\qquad k=1,\dots,N
~,
\label{bkt3}
\ee
where the rescaled frequency targets have coordinates
${s_k''}^{(i)} := h_i {s_k'}^{(i)} = h_i\gamma_i s_k^{(i)}$, $i=1,\dots,d$.
Intuitively, the factor $h_i$ arises because the fine grid of spacing $h_i$
has to be stretched to unit spacing to be interpreted as a Fourier series.

{\bf Step 3 (correction).}
Finally, as in type 1, in order to compensate for the spreading in step 1
(in primed coordinates) a diagonal correction is needed,
$$
\tilde f_k = p_k b_k,\quad 
p_k = h_1\dots h_d \hat\psi(\sss'_k)^{-1} = (2/w)^d
\bigl(
\hat\phi(\alpha_1{s'_k}^{(1)}) \cdots \hat\phi(\alpha_d{s'_k}^{(d)})
\bigr)^{-1}
\!,
\quad k=1,\dots,N.
$$
But, in contrast to types 1 and 2,
the set of frequencies at which $\hat\phi$ must be evaluated is 
{\em \NU}, so there is no phase winding trick.
Rather, $dpN$ cosines must be evaluated,
recalling that $p$ is the the number of positive quadrature nodes.
Despite this cost, this step consumes only a small fraction of the
total computation time.

\begin{rmk}
  Empirically, we find that
  using the same overall requested precision $\eps$
  in the above steps 1 and 2 gives overall error still close to $\eps$.
  It has been shown in 1D (see term $E_3$ in \cite[p.~45]{elbel}) that
  the type 3 error is bounded by the error in performing the above step 2
  multiplied by $\rmax$,
  the dynamic range of $\hat\psi$ over the usable frequency band
  (see Fig.~\ref{f:spreadalias}(b)).
  Using $n\approx \rat N$,
  \eqref{EShat1}, and \eqref{gam} with $\gamma\approx1$ we approximate
  \be
  \rmax := \frac{\hat\psi(0)}{\hat\psi(N/2)}
  = \frac{\hat\phi(0)}{\hat\phi(\pi w/2\rat)}
  \approx e^{\freq - \sqrt{\freq^2 - (\pi w/2\rat)^2}}
  = e^{\bigl(1 - \sqrt{1-(2\rat-1)^{-2}}\bigr)\freq}
  ~,
  \label{rmax}
  \ee
  which for $\rat=2$ gives $\rmax \approx e^{0.057\,\freq}$.
  From \eqref{params}, $\beta\le36$ for any $\eps\ge10^{-15}$, so $\rmax\le 8$,
  which is quite small. This helps to justify the above choice of tolerances.
\end{rmk}

{\bf Choice of fine grid size (``Step 0'' for type 3).}
Finally we are able to give the recipe for choosing the fine grid
sizes $n_i$ (which of course in practice precedes the above three steps).
This relies on aliasing error estimates \cite{elbel}
for steps 1 and 3 that we explain here only heuristically.
In section~\ref{s:anal}
we will see that spreading onto a uniform grid
of size $h_i$ induces a lattice of
aliasing images separated by $n_i$ in frequency space,
so that the correction step is only accurate to precision $\eps$ out to
frequency magnitude $n_i/2\rat$.
Thus, since $|{\sss'_k}^{(i)}| \le \gamma_i S_i$ for all $i$ and $k$,
the condition
\be
\gamma_i S_i \le \frac{n_i}{2\rat}
~,\qquad i=1,\dots,d
\label{cond2}
\ee
is sufficient.
Combining \eqref{cond1} and \eqref{cond2}, then solving as equalities
for the smallest $n_i$
gives the recipe for the optimal parameters (similar to \cite[Rmk.~1]{nufft3}),
\be
n_i \;=\; \frac{2\rat}{\pi}X_iS_i + w
~,\qquad
\gamma_i \;=\; \frac{n_i}{2\rat S_i}
~,
\qquad i=1,\dots,d
~.
\label{ng}
\ee

\begin{rmk}[FFT size for type 3]
  The product of the grid sizes $n_i$ in each dimension $i=1,\dots,d$
  sets the number of modes, hence the FFT effort required,
  in the type 2 transform in step 2.
Crucially,  this is independent of the numbers of sources $M$ and of
  targets $N$.
  Rather, $n_i$ scales like the space-frequency product $X_iS_i$.
  This connects to the Fourier uncertainty principle:
  $n_i$ scales as the number of ``Heisenberg boxes''
  needed to fill the centered rectangle enclosing the data.
  In fact, since the number of degrees of freedom \cite[p.~391]{slepianrev}
  (or ``semiclassical basis size'' \cite{davisheller})
  needed to represent functions
  living in the rectangle $[-X_i,X_i]\times[-S_i,S_i]$ is its area divided
  by $2\pi$, namely $2X_iS_i/\pi$, we see that
  $n_i$ is asymptotically $\rat$ times this basis size.
  \label{r:heis}
\end{rmk}

{\bf Efficiently handling poorly-centered data.}
The above remark shows that
the type 3 is helped by translating all coordinates
$\xx_j$ and $\sss_k$ so that their respective bounding boxes are centered around
the origin. This reduces the bounds $X_i$ and $S_i$ defined
by \eqref{XS}, hence reduces $n_i$ and thus the cost of the FFT.
Translations in $\xx$ or in $\sss$ are cheap to apply using the factorization
\be
\sum_{j=1}^M c_j e^{i (\sss_k+\sss_0) \cdot (\xx_j+\xx_0)}
\;= \;
e^{i(\sss_k+\sss_0)\cdot\xx_0} \sum_{j=1}^M (e^{i\sss_0\cdot \xx_j}c_j) e^{i \sss_k\cdot\xx_j}~.
\label{trans}
\ee
Thus the type 3 transform for translated data
can be applied by pre-phasing the strengths by $e^{i\sss_0\cdot \xx_j}$,
doing the transform, then post-multiplying by $e^{i(\sss_k+\sss_0)\cdot\xx_0}$.
The extra cost is $\bigO(N+M)$ complex exponentials.
In our library, if input or output points are sufficiently poorly centered,
we apply \eqref{trans} using as $\xx_0$ or $\sss_0$ the
means of the minumum and maximum
coordinates in each dimension.

\begin{rmk}[type 3 efficiency]
  Remark~\ref{r:heis} also shows that input data can be chosen for which the
  algorithm is arbitrarily inefficient.
  For example, with only two points ($M=N=2$) in 1D with
  $x_1=-X$, $x_2=X$, $s_1=-S$, $s_2=S$,
  then by choosing $XS$ huge, \eqref{ng} implies that the algorithm
  will require a huge amount of memory and time.
  Obviously in such cases a direct summation of \eqref{3} is preferable.
  However, for $N$ and $M$ large but with clustered data,
  a butterfly-type algorithm which hierarchically exploits \eqref{trans}
  could be designed; we leave this for future work.
\end{rmk}

\section{Error analysis and parameter choices}
\label{s:err}

Here we summarize a rigorous estimate (proven in \cite{nufftanal})
on the aliasing error of the 1D type 1 and 2 algorithms of
section~\ref{s:alg}, when performed in exact arithmetic.
We then use this to justify the algorithm parameter choices
stated in \eqref{params}.
Finally we evaluate and discuss the gap between this estimate and empirical
errors.

\subsection{Theoretical results for the ES kernel}
\label{s:anal}

Let $\mbf{f}$ be the vector of $f_k$ outputs defined by \eqref{1} in 1D,
and $\tilde{\mbf{f}}$ be the analogous
output of the above type 1 NUFFT algorithm
in exact arithmetic. We use similar notation for type 2.
By linearity, and the fact that the type 2 algorithm is the adjoint of the
type 1, the output aliasing error vectors must take the form
\be
\tilde{\mbf{f}} - \mbf{f} = E\mbf{c}
\quad \mbox{ (type 1)~, } \qquad\qquad
\tilde{\mbf{c}} - \mbf{c} = E^\ast\mbf{f}
\quad \mbox{ (type 2)~, }
\label{E}
\ee
for some matrix $E$.
Standard analysis
\cite[(1.16)]{nfftchap}
\cite[(4.1)]{fourmont} \cite[Sec.~V.B]{fessler}
(or see \cite{nufftanal})
involving the Poisson summation formula
shows that, with the choice \eqref{pk1} for $p_k$,
$E$ has elements
\be
E_{kj} = g_k(x_j) ~,
\qquad \mbox{where }
\quad
g_k(x) = \frac{1}{\hat\psi(k)}\sum_{m\neq 0} \hat\psi(k+mn) e^{i(k+mn)x}
~.
\label{gkx}
\ee
Since $|k|\le N/2$,
error is thus controlled by a phased sum of the tails of $\hat\psi$
at frequency magnitudes at least $n-N/2$; see Fig.~\ref{f:spreadalias}(b).

It is usual in the literature to seek a uniform bound on elements of $E$
by discarding the information about $x_j$, so that $|E_{kj}|\le\eps_\infty$,
$\forall kj$, where
\be
\eps_\infty \;:=\; \max_{|k|\le N/2}
\|g_k\|_\infty \;\le\;
\frac{
  \max_{|k|\le N/2, \,x\in\RR} \left|
    \sum_{m\neq 0} \hat\psi(k+mn) e^{i(k+mn)x}
  \right|
}{\min_{|k|\le N/2} |\hat\psi(k)|}
~.
\label{epsest}
\ee
The latter inequality is close to tight because
$\rmax$ defined by \eqref{rmax}
controls the loss due to bounding numerator and denominator separately,
and is in practice small.
\begin{rmk}
A practical heuristic for $\eps_\infty$
is sketched in Fig.~\ref{f:spreadalias}(b):
assuming that i) $\hat\psi(k)$ decreases monotonically
with $|k|$ for $|k|\le N/2$, and  ii) the worst-case sum (numerator in \eqref{epsest}) is dominated by the single value with smallest $|k|$,
then we get $ \eps_\infty \approx |\hat\psi(n-N/2) / \hat\psi(N/2) |$,
whose logarithm is shown in the figure.
\label{r:heur}
\end{rmk}    

A common use for \eqref{epsest}
is the simple $\ell_1$-$\ell_\infty$ bounds for \eqref{E}
(see \cite{steidl98} \cite[p.12]{fessler}):
\be
\max_{k\in\KK}|\tilde f_k - f_k| \le \eps_\infty \|\mbf{c}\|_1
\quad \mbox{ (type 1), } \quad
\max_{1\le j\le M}|\tilde c_j - c_j| \le \eps_\infty \|\mbf{f}\|_1
\quad \mbox{ (type 2). }
\label{1nrm}
\ee


These results apply to any spreading kernel; we now specialize to the
ES kernel.
Fix an upsampling factor $\rat>1$.
Given a kernel width $w$ in sample points, one must choose in \eqref{ES}
an ES kernel parameter $\freq$ such that $\hat\psi$ defined in \eqref{psi1}
has decayed to its
exponentially small region once the smallest aliased
frequency $n-N/2 = n(1-1/2\rat)$ is reached; see \eqref{epsest} and
Fig.~\ref{f:spreadalias}(b).
To this end we fix a ``safety factor'' $\gamma$, and set
\be
\freq \;=\; \freq(w) \;:=\; \gamma \pi w (1-1/2\rat)~,
\label{gam}
\ee
so that for $\gamma=1$ the exponential cutoff occurs exactly at $n-N/2$,
while for $\gamma<1$ the cutoff is safely smaller than $n-N/2$.
With this set-up, the following states that
the aliasing error converges almost exponentially
with respect to the kernel width $w$.

\begin{thm}[\cite{nufftanal}] 
  For the 1D type 1 and 2 NUFFT, fix $N$ and $\rat$ (hence the upsampled grid $n=\sigma N$) and the safety factor $\gamma\in(0,1)$.
  With $\freq(w)$ as in \eqref{gam},
  then the aliasing
  error uniform bound \eqref{epsest} converges with respect to
  kernel width $w$ as
  \be
  \eps_\infty \;=\; \bigO\left( \sqrt{w} e^{-\pi w \gamma
    \sqrt{1-1/\rat - (\gamma^{-2}-1)/4\rat^2}}
  \right)
  ~, \qquad w\to \infty ~.
  \label{ESerr}
  \ee
  \label{t:ESerr}
\end{thm}  
Its somewhat involved proof is detailed in \cite{nufftanal}.
A key ingredient is that,
asymptotically as $\freq\to\infty$,
$\hat\phi$ has the same ``exponential of semicircle'' form
(up to algebraic factors)
in the below-cutoff domain
$(-\freq,\freq)$ that $\phi$ itself has in $(-1,1)$;
compare Figs.~\ref{f:kernel}(b) and (c).
Specifically, fixing the scaled frequency $\rho\in(-1,1)$,
\cite{nufftanal} proves that,
\be
\hat\phi(\rho\freq) \; =\;
\sqrt{\frac{2\pi}{\beta}} \frac{1}{(1-\rho^2)^{3/4}} e^{\freq(\sqrt{1-\rho^2}-1)}
\left[ 1 + \bigO(\freq^{-1}) \right]
~, \qquad \freq\to\infty
~.
\label{EShat1}
\ee

\begin{rmk}[Comparison to \KB\ bounds]
  \label{r:fourmont}
  In the limit $\gamma\to1^{-}$, \eqref{ESerr}
  approaches the same exponential rate as \eqref{KBerr},
  and with an algebraic prefactor improved by a factor $\sqrt{w}$.
  On the other hand, \eqref{KBerr} has an explicit constant.
\end{rmk}

\bfi 
\mbox{
  \ig{width=1.72in}{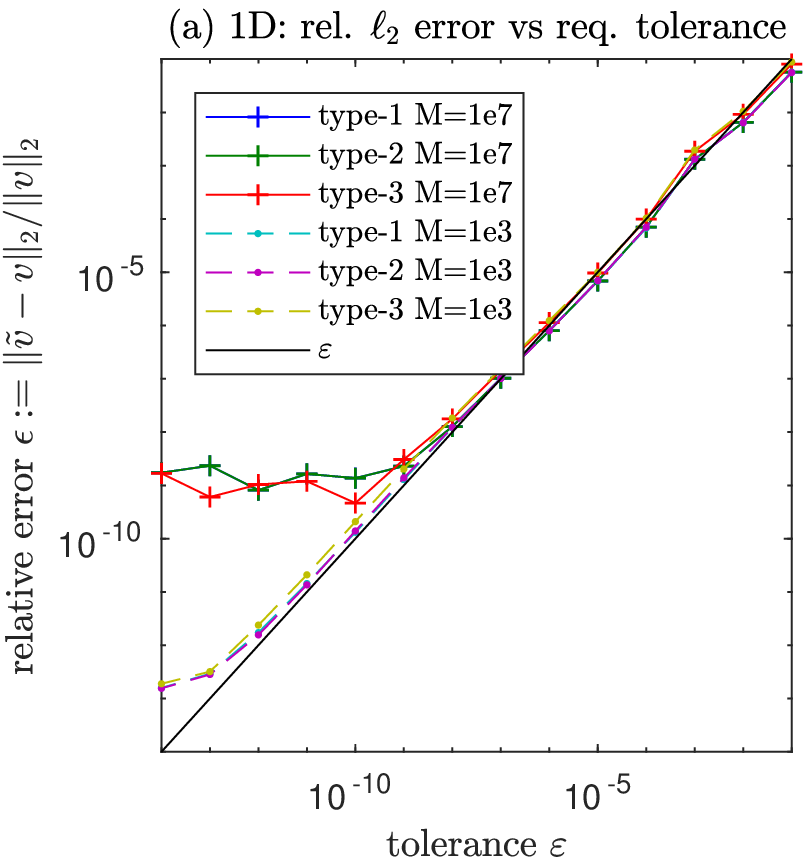}
  \ig{width=1.67in}{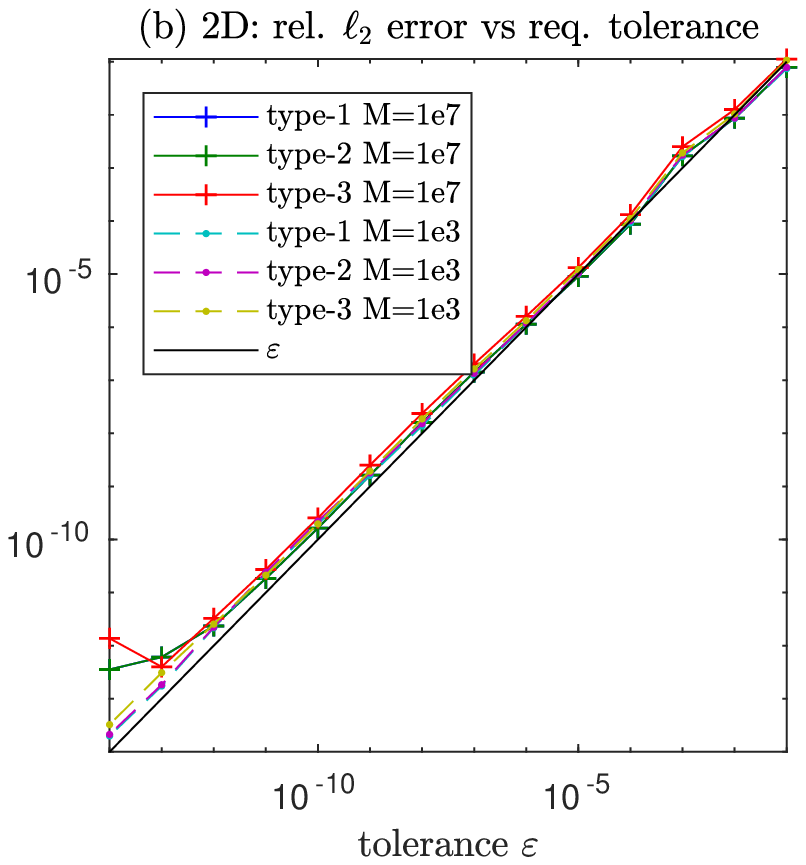}
  \ig{width=1.67in}{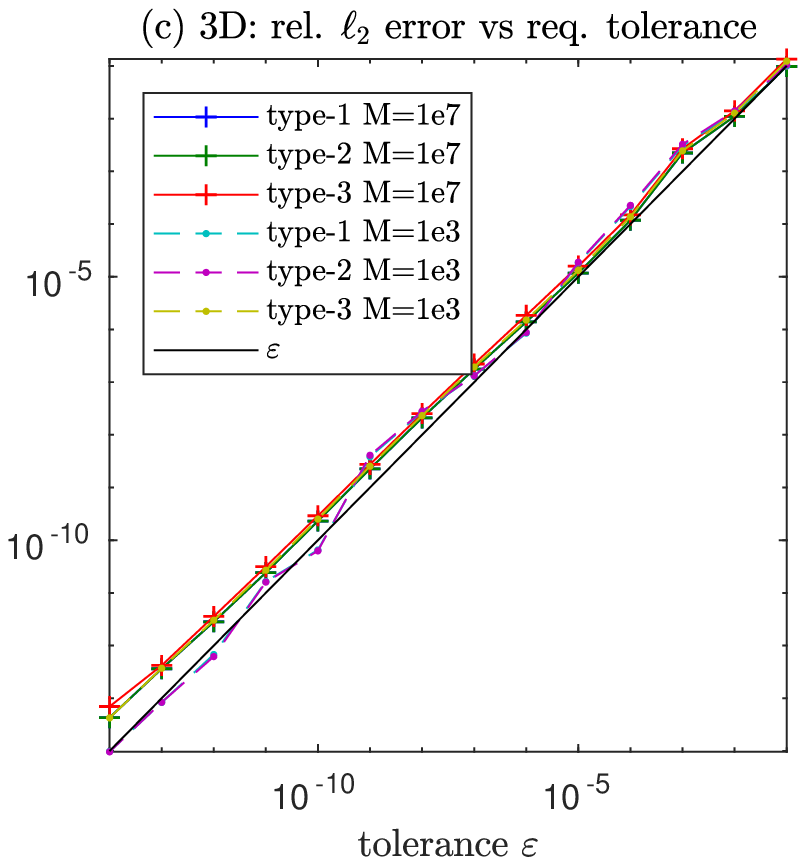}
  }
\ca{Comparison of empirical relative $\ell_2$ error ($\epsilon$) vs
  the requested tolerance ($\eps$), for all nine FINUFFT routines.
  In each case $N\approx M$, with $N_i$ roughly equal,
  with two problem sizes included ($M=10^3$ and $10^7$).
  Nonuniform points are uniformly randomly distributed in $[0,2\pi)^d$,
    while strength data is complex Gaussian random.
  $v$ denotes the true output vector, either $\mbf{f}$ or $\mbf{c}$
  (this is computed with tolerance $10^{-15}$), and $\tilde v$ the computed
  output vector at requested tolerance $\eps$.
  }{f:acc}
\efi

\subsection{ES kernel parameter choices and empirical error}
\label{s:params}

We now justify and test the parameter choices \eqref{params}.
With $\rat=2$, the factor 2.30 in \eqref{params}
corresponds to a satefy factor $\gamma\approx0.976$ in \eqref{gam}, very close to 1.
Note that $\gamma=1$ would give a factor $2.356$; however,
we find that the factor $2.30$ gives a slightly lower typical error for a given $w$
than pushing $\gamma$ closer or equal to 1
(this is likely due to the continued drop for $\xi>\beta$ visible in
Fig.~\ref{f:kernel}(c)).
Fessler \etal found a similar factor 2.34 when optimizing the KB kernel
\cite[Fig.~11]{fessler}.

The width $w$ is set by solving \eqref{ESerr} with its algebraic prefactor
dropped, to give
$w \approx
|\log \eps_\infty|/\pi\gamma\sqrt{1 - 1/\rat - (\gamma^{-2}-1)/4\rat^2}
\,+\,\mbox{const}$.
Interpreting $\eps_\infty$ as the requested tolerance $\eps$ (see
Remark~\ref{r:err} below),
and inserting $\gamma\approx0.976$, gives
\be
w \;\approx\; 1.065 |\log_{10} \eps| \; + \mbox{ const}
~.
\label{wslope}
\ee
As we show below, the factor 1.065 may be replaced by unity
while still giving empirical errors close to the requested tolerance.
The constant term in \eqref{wslope} is fit empirically.
Thus we have explained the parameter choices \eqref{params}.

In many applications one cares about relative $\ell_2$ error in the output
vector, which we will denote by
\be
\epsilon := \frac{\|\tilde{\mbf{f}}-\mbf{f}\|_2}{\|\mbf{f}\|_2}
\quad \mbox{ (types 1 and 3)~, } \qquad
\epsilon := \frac{\|\tilde{\mbf{c}}-\mbf{c}\|_2}{\|\mbf{c}\|_2}
\quad \mbox{ (type 2)~. }
\ee
Following many \cite{dutt,nfftchap}, we will use this metric for testing.
Fig.~\ref{f:acc} measures this metric
for FINUFFT for all nine transform types at two different problem sizes,
with random data and randomly located \NU\ points.
This shows that, using the choice \eqref{params},
the achieved relative error $\epsilon$ well matches the requested tolerance
$\eps$, apart from when round-off error dominates.
The mean slope of the logarithm of the empirical error $\epsilon$ with
respect to that of $\eps$ in Fig.~\ref{f:acc} is slightly less than
unity, due to the design choice of approximating the slope 1.065 in \eqref{wslope} by unity in \eqref{params}.

\begin{rmk}[Rounding error]
  Double-precision accuracy is used for all machine calculations in
  the library by default, and also in the studies in this work.
  The resulting rounding error is only apparent above aliasing error
  for the large 1D and 2D transforms at high accuracy.
  Fig.~\ref{f:acc}(a) shows that our library's
  accuracy is limited to 9 relative digits
  in 1D for $M=N=10^7$; more generally
  Fig.~\ref{f:1d} shows that rounding error is similar, and essentially
  the same for all tested libraries.
  When $M\approx N$ we find, in 1D, that the rounding
  contribution to $\epsilon$ is roughly
  $N$ times machine precision.
  Taking into account their choice of error norm,
  this concurs with the findings of \cite[Fig.~2.3]{townsendnufft}.
  See \cite[\S 1.4]{nfftchap} for NUFFT rounding error analysis in 1D.
  We observe in 2D and 3D that
  it is $\max_{i} N_i$ that scales the rounding error;
  thus, as Fig.~\ref{f:acc}(b-c) shows, it is largely irrelevant even
  for large $M$.
  \label{r:rounding}
\end{rmk}

Finally, we discuss the gap between
any theoretical aliasing error estimate deriving from
\eqref{epsest}---%
this includes
\eqref{KBerr} and \eqref{ESerr}---%
and the {\em empirical} relative $\ell_2$ error $\epsilon$ due to aliasing.
The best possible type 1 bound from \eqref{E}--\eqref{epsest}
is via the Frobenius norm $\|E\|_F \le \sqrt{MN}\eps_\infty$,
so
$$
\|\tilde{\mbf{f}}-\mbf{f}\|_2 \;\le\; \sqrt{MN}\eps_\infty \|\mbf{c}\|_2
$$
Writing the transform \eqref{1} as $\mbf{f}=A\mbf{c}$, where
$A$ has elements $A_{kj} =e^{ikx_j}$,
this gives
\be
\epsilon \;\le\; \sqrt{MN}\eps_\infty \frac{\|\mbf{c}\|_2}{\|\mbf{f}\|_2}
\;\le\; \sqrt{MN}\eps_\infty \frac{1}{\sigma_\tbox{min}(A)}
\label{worst}
\ee
where in the last step the best bound applying to all nontrivial $\mbf{c}$
is used, and $\sigma_\tbox{min}(A)$ denotes the smallest singular value of $A$,
or zero if $M>N$.
Thus if $M>N$, there cannot be a general type 1 bound on $\epsilon$, simply
because, unlike the DFT, the {\em relative condition number} of the type 1 task \eqref{1}
may be infinite (consider $\mbf{c} \in \mbox{Nul }A$,
so $\mbf{f}=\mbf{0}$).%
\footnote{The condition number may also be huge even if $M\le N$.
  The following MATLAB code,
  in which \NU\ points lie randomly in only {\em half} of
  the periodic interval,
  outputs typically $10^{-15}$:
  \quad {\tt M=80; N=100; A = exp(1i*(-N/2:N/2-1)'*pi*rand(1,M)); min(svd(A))}
}

However there are two distinct mechanisms by which \eqref{worst} is pessimistic
in real-world applications:
\ben
\item
  For {\em typical} input data, $\|\mbf{f}\|_2$ is not smaller
  than $\|\mbf{c}\|_2$; in fact
  (as would be expected from randomized phases in $A$)
  typically $\|\mbf{f}\|_2 \approx \sqrt{N} \|\mbf{c}\|_2$.
  The growth factor is close to $\sigma_\tbox{max}(A)$. Thus
  the problem is generally well-conditioned.
  See B\"ottcher--Potts \cite[Sec.~4]{boettcher} for a formalization
  in terms of ``probabilistic condition number''.
\item The uniform bound \eqref{epsest}
  discards phase information
  in the elements of $E$, which in practice induce large cancellations
  to give errors that are improved by a factor $\sqrt{M}$
  over bounds such as \eqref{1nrm}.
  Such ideas enable improved aliasing error estimates
  in a looser norm:
  \eg interpreting \eqref{2} as point samples of a {\em function}
  $c(\xx) = \sum_{\kk\in\KK} f_\kk e^{-i\kk\cdot\xx}$,
  the error of the latter is easily bounded in $L^2([-\pi,\pi]^d)$,
  as clarified by Nestler~\cite[Lemma~1]{nestler}
  (also see \cite{l2jacob}).
\een

\begin{rmk}[Empirical $\ell_2$ relative error]
  Assuming both mechanisms above apply in practice,
  they can be combined to replace \eqref{worst} with the heuristic
  $$\epsilon \;\approx \;\frac{\sqrt{N}\eps_\infty \|\mbf{c}\|_2}{\sqrt{N}\|\mbf{c}\|_2} \;= \; \eps_\infty~,$$
  which justifies the interpretation of the requested tolerance $\eps$ 
  as the uniform bound $\eps_\infty$ in \eqref{epsest} when setting kernel parameters.
  The result is that, as in Fig.~\ref{f:acc}, barring rounding error,
  relative error $\epsilon$ is almost always similar to the
  tolerance $\eps$.
  \label{r:err}
\end{rmk}

To summarize,
rather than design the kernel parameters around the rigorous
but highly pessimistic worst-case analysis \eqref{worst},
we (as others do) design for typical errors.
Thus, before trusting the relative error,
the user is recommended to
check that for their input data
the desired convergence with respect to $\eps$ has been achieved.

\section{Implementation issues}
\label{s:imp}

Here we describe the main software aspects that accelerated the library.
The chief computational costs in any NUFFT call
are the spreading (types 1 and 3) or interpolation (type 2),
scaling as $\bigO(w^dM)$, and the FFT, scaling as $\bigO(N \log N)$.
We use the multithreaded FFTW library
for the latter, thus in this analysis we focus on spreading/interpolation.
In comparison, the correction steps, as explained in section~\ref{s:alg},
are cheap: 1D evaluations of $\hat\psi$
are easily parallelized over the $p$ quadrature nodes
(or, for type 3, the frequency points) with OpenMP, and the correction
and reshuffling of coefficients is memory-bound so does not
benefit from parallelization.

\subsection{Bin sorting of \NU\ points for spreading/interpolation}

When $N$ is large, the upsampled grid (with $\rat^dN$ elements) is
too large to fit in cache.
Unordered reads/writes to RAM are very slow (hundreds of clock cycles),
thus looping through the \NU\ points in an order which preserves
locality in RAM uses cache well and speeds up spreading and interpolation,
by a factor of typically 2--10, including the time to sort.%
\footnote{This is illustrated
  by running {\tt test/spreadtestnd 3 1e7 1e7 1e-12 x 0 1} where {\tt x} is
  {\tt 0} (no sort) or {\tt 1} (sort). For interpolation ({\tt dir=2}),
  we find a speedup factor 6 on a Xeon, or 14 on an i7.}
Each \NU\ point requires accessing a block extending $\pm w/2$ grid points in
each dimension, so there is no need to sort to the nearest grid point.
Thus we set up boxes of size 16 grid points
in the fast ($x$) dimension, and, in 2D or 3D,
size 4 in the slower ($y$ and $z$) dimensions. These sizes are
a compromise between empirical speed and additional RAM needed for the sort.
Then we do an ``incomplete'' histogram sort:
we first count the number of points
in each box and use this to construct the breakpoints between bins,
then write point indices lying in each box into that bin,
finally reading off the indices in the box ordering (without sorting
inside each bin).
This bin sort is multi-threaded, except for the low-density case $M<N/10$
where we find that the single-thread version is usually faster.

\bfi 
\centering
\ig{width=5.2in}{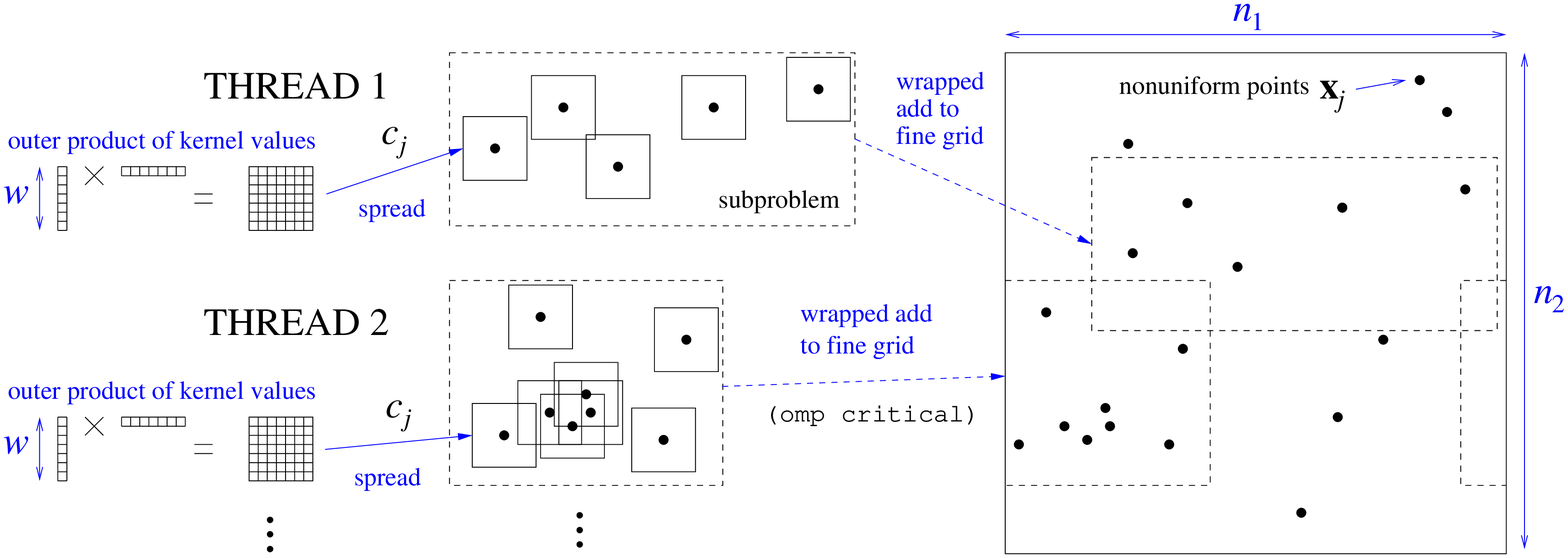}
\ca{Sketch of parallel load balanced spreading scheme used for
  type 1 and type 3 transforms, showing the 2D case. Only
  two of the threads are shown.
}{f:t1blocks}
\efi

\subsection{Parallel spreading}
\label{s:block}

The interpolation task \eqref{interp} parallelizes well with OpenMP,
even for highly nonuniform distributions.
Each thread is assigned a subset of the points $\{\xx_j\}$;
for each point it reads a block of size $w^d$ from the fine grid,
does a weighted sum using as weights the tensor product of 1D kernels,
then writes the sum to a distinct output $c_j$.

In contrast, spreading \eqref{bl} 
{\em adds} to blocks in the fine grid.
Blocks being overwritten by different threads may collide,
and, even if atomic operations are used to ensure a correct result,
{\em false sharing} of cache lines \cite[Sec.~5.5.2]{usingopenmp}
makes it inefficient.
A conventional solution is to assign threads equal distinct slices of the fine
grid \cite{nfft}; however, for nonuniform distributions
this could lead to arbitrarily poor load balancing.
Instead, we group sorted points $\{\xx_j\}$ into ``subproblems'' of size up to
$10^4$, which are assigned to threads; see Fig.~\ref{f:t1blocks}.
(This choice is heuristic;
an optimal size would depend on L3 cache and the number of threads.)
To handle a subproblem, a thread finds the cuboid bounding all
$\{\xx_j\}$ in the subproblem, allocates a local cuboid of this size,
spreads onto the cuboid, and finally adds the cuboid back into
the fine grid.
Since subproblems may overlap on the latter grid,
this last step needs a {\tt \#pragma omp critical} block
to avoid collisions between writes; however this causes minimal
overhead since almost all the time is spent spreading to cuboids.
The scheme is adaptive:
regardless of the point distribution, all threads are kept busy almost
all of the time.
The scheme requires additional RAM of order the fine grid size.

A further advantage is that the periodic wrapping of grid indices,
which is slow, may be avoided: cuboids are padded by $w/2$ in each
dimension and written to without wrapping. Index wrapping is only
used when adding to the fine grid.



\subsection{Piecewise polynomial kernel approximation}
\label{s:poly}

The 1D ES kernel \eqref{ES} requires one real-valued {\tt exp}
and {\tt sqrt} per evaluation.
However, we find that the throughput 
depends drastically on the CPU type (i7 or Xeon), compiler
(GCC vs Intel ICC), and kernel width $w$ (the inner loop length
that the compiler may be able to vectorize via SIMD instructions).
For instance, a Xeon E5-2643 (with AVX2) with GCC version $\le 7$.x
achieves only 40M evals/sec/thread,
while the same CPU with ICC gives 50--200M evals/sec/thread.
We believe this is due to 
compiler-provided {\tt exp} instructions that exploit SIMD.
Similar variations occur for the i7.
Seeking a reliably efficient kernel evaluation on open-source compilers
(\eg GCC), we replaced the kernel evaluation by a polynomial look-up table.
The result gives 350--600M evals/sec/thread,
and accelerates 1D spreading/interpolation by a factor 2--3
(the effect in $d=2,3$ is less dramatic).

The look-up table works as follows. For each \NU\ point coordinate,
the 1D kernel $\psi(x)$ must be evaluated at $w$ ordinates
$x, x+h, \ldots,x+(w-1)h$
(see Fig.~\ref{f:spreadalias}(a)).
We break the function $\psi$ into $w$ equal-width intervals and approximate
each by a centered polynomial of degree $p$.
Ordinates within each of the $w$ intervals are then the same, allowing
for SIMD vectorization. 
The approximation error need only be small relative to $\eps_\infty$;
we find $p=w+3$ suffices.
Monomial coefficients are found by solving a Vandermonde system on
collocation points on the boundary of a square in the complex plane
tightly enclosing the interval.
For each of the two available upsampling factors ($\rat=2$ and $5/4$),
and all relevant $w$, we automatically
generate C code containing all coefficients and Horner's evaluation scheme.

\begin{rmk}
  Piecewise polynomial approximation could confer on any kernel
  (\eg KB, PSWF) this same high evaluation speed.
  However, AVX512 and future SIMD instruction sets may accelerate
  {\tt exp} evaluations, making ES even faster.
  Since we do not know which will win with future CPUs,
  our library uses the ES kernel.
\end{rmk}

\bfi  
\mbox{
\ig{height=1.7in}{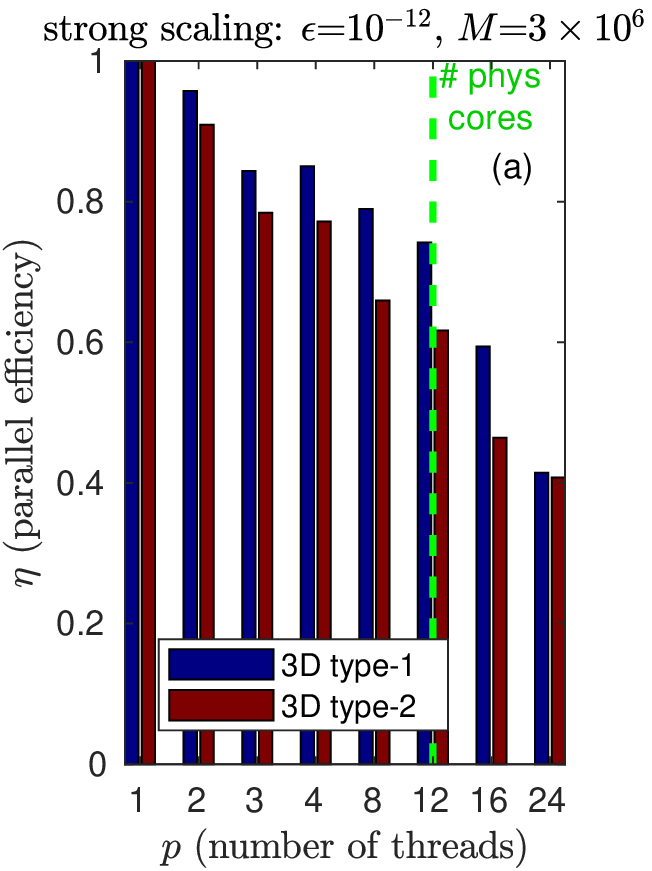}
\ig{height=1.7in}{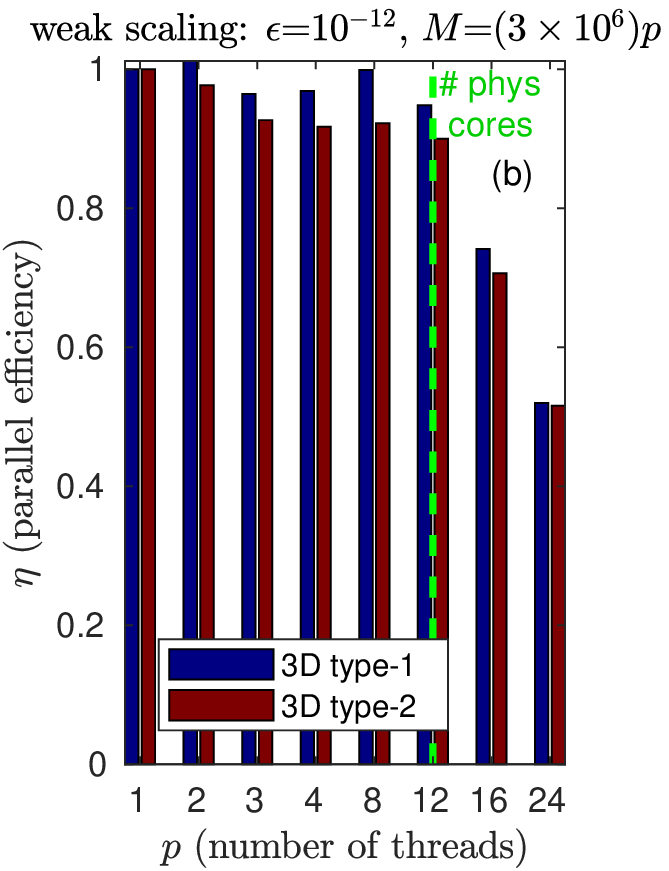}
\ig{height=1.7in}{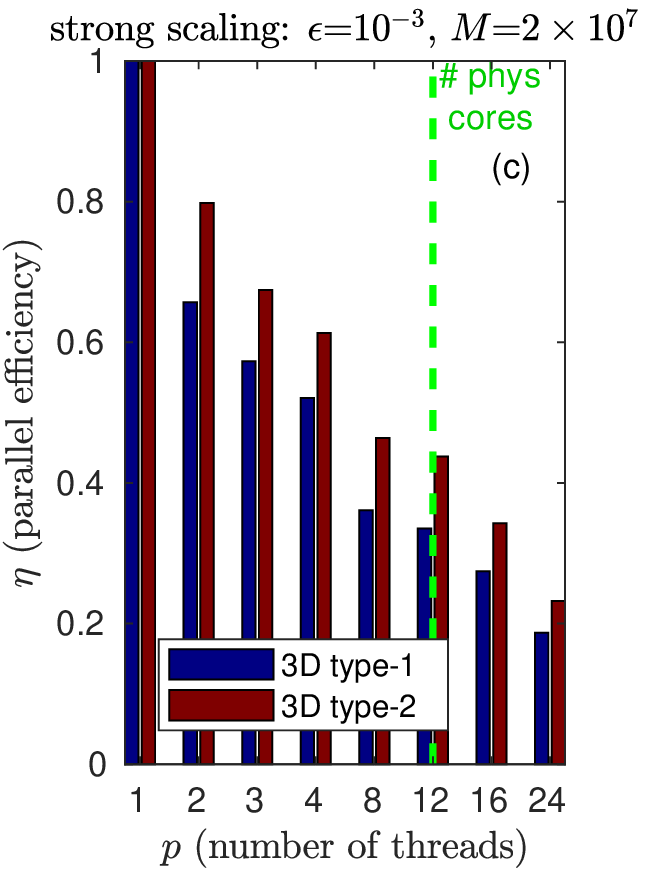}
\ig{height=1.7in}{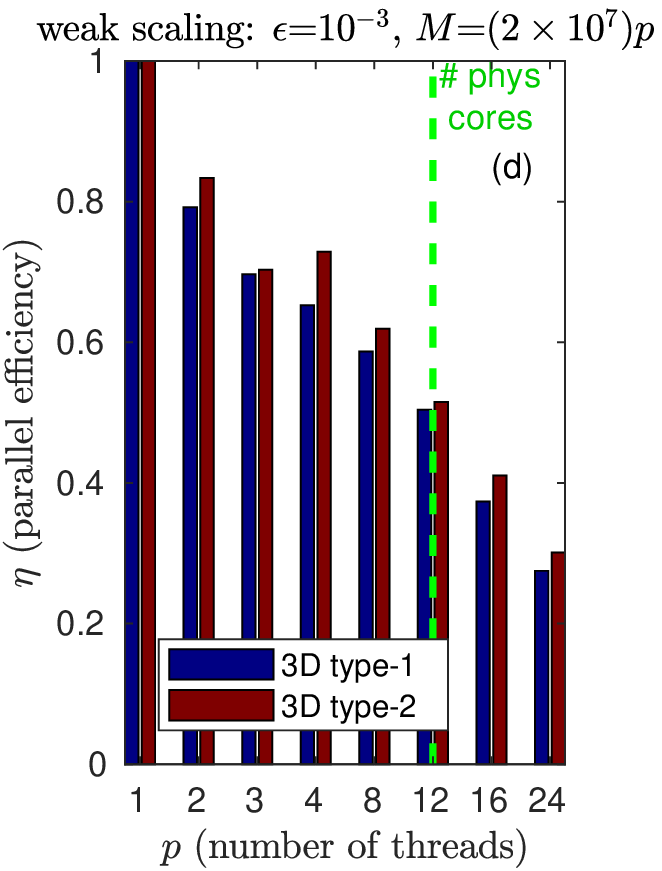}
}
\ca{Parallel scaling of FINUFFT in 3D, for $p$ threads of
  a Xeon desktop with 12 physical cores.
  Both type 1 and type 2 tasks are tested.
  In all cases there are $N=100^3$ Fourier modes, and $M$,
  the number of ``sph quad'' nodes (see section~\ref{s:perf}), is as shown.
  (a--b) show a high-accuracy case (12 digits).
  (c--d) show low accuracy (3 digits).
  For strong scaling the efficiency is the speed-up factor divided by $p$;
  for weak it is the speed-up factor
  for a problem size $M$ proportional to $p$.
}{f:parsca}
\efi

\section{Performance tests}
\label{s:perf}

Tests were run on a desktop with two Intel Xeon 3.4 GHz E5-2643
CPUs (each with 20 MB L3 cache),
giving 12 total physical cores (up to 24 threads), and 128 GB RAM, running EL7
linux.
Unless specified we compile all codes with
GCC v.7.3.0.
We compiled FINUFFT version 1.0 with flags
{\tt -fPIC -Ofast -funroll-loops -march=native}.
Experiments were driven using the MEX interface to MATLAB R2016b.
In the codes that use FFTW (\ie FINUFFT and NFFT), we use version 3.3.3
and set its plan method to
\verb+FFTW_MEASURE+ (see \cite{fftw}),
which sometimes takes a very long time to plan during the first call.
Thus, to show realistic throughput we time only subsequent calls,
for which FFTW looks up its stored plan.
To minimize variation we take the best of the three subsequent calls.

{\bf Tasks.}\;
To assess the efficiency of our contributions---rather than merely
measure the speed of FFTW---we choose ``high density'' tasks where $M$ is
somewhat larger than $N$, so the
FFT is at most a small fraction of the total time.
In 1D, since timings do not vary much with point distribution,
we always test with $x_j$ iid uniform random in $[-\pi,\pi]$.
For $d=2,3$ we use the following distributions
(see insets in Figs.~\ref{f:2d}--\ref{f:3dmulti}):
\bi
\item 2D ``disc quad'': a polar grid over the disc of radius $\pi$,
  using roughly $\sqrt{M}$ radial Gauss--Legendre nodes and
  $\sqrt{M}$ equi-spaced angular nodes.
\item 3D ``rand cube'': iid uniform random in $[-\pi,\pi]^3$.
\item 3D ``sph quad'': a spherical grid in the ball of radius $\pi$,
  using $\sqrt{M}/2$ radial Gauss--Legendre nodes and a
  $\sqrt{M} \times 2\sqrt{M}$ tensor-product grid on each sphere.
\ei
Here the first and last are realistic quadrature schemes for NUFFT
applications \cite{cryo}. They involve a divergence in point density
at the origin of the form $r^{1/2-d}$ for $d=2,3$.
We choose input strengths or coefficients as
iid complex Gaussian random numbers.

{\bf Parallel scaling.}\;
Fig.~\ref{f:parsca} shows parallel scaling tests of 3D type 1 and 2 FINUFFT.
The highly-nonuniform ``sph quad'' distribution was used
in order to test the load balancing described in section~\ref{s:block}.
For $\eps=10^{-12}$, where each point interacts with $13^3=2197$ fine
grid points, weak scaling
(where $M$ grows with $p$ the number of threads)
shows 90\% parallel efficiency for $p\le 12$ (one thread per physical core).
Above this, hyper-threading is used:
as expected, although it provides a slight net speed boost,
measured in threads its
parallel efficiency falls far short of that for $p\le 12$.
Strong scaling (acceleration at fixed $M$) is a tougher test,
dropping to 62--74\% at $p=12$.

For a low-accuracy test ($\eps=10^{-3}$), the kernel is narrower, touching only
$4^3=64$ fine grid points, thus the RAM access pattern is more
random relative to the number of FLOPs.
We believe the resulting lower parallel efficiencies are due to
memory bandwidth, rather than flops, being the bottleneck.
That said, at $p=24$ threads (full hyper-threading), both transforms
are still 5--7 times faster than for a single core.

\begin{table} 
\centering
  {\footnotesize 
  \begin{tabular}{l|llllll}
Code name & kernel & langauge & on-the-fly & omp & periodic domain & notes
\\
\hline
FINUFFT & ES & C++ & yes & yes & $[-\pi,\pi]^d$ & 
\\  
CMCL & Gaussian & Fortran & yes & no & $[-\pi,\pi]^d$ &
\\
NFFT & bkwd.\ KB & C & yes or no & yes & $[-1/2,1/2]^d$ &
\\
MIRT & optim.\ KB & MATLAB & no & no & $[-\pi,\pi]^d$ &
\\
BART & KB, $w=3$ & C & yes & yes & $\prod_{i=1}^d [-N_i/2,N_i/2]$ & $d=3$ only
\\
\hline
  \end{tabular}
}
\vspace{1ex}
\ca{Summary of NUFFT libraries tested.
  ``kernel'' lists the default spreading function $\psi$
  (some allow other kernels). ``language'' is the coding language.
  A ``yes'' in the column ``on-the-fly'' indicates that no precomputation/plan
  phase is needed, hence low RAM use per nonuniform point.
  ``omp'' shows if multithreading (e.g.\ via OpenMP)
  is available.
  See sections~\ref{s:prior} and \ref{s:bench}
  for more details.}{t:codes}
\end{table}   

\bfi
\ig{height=2.2in}{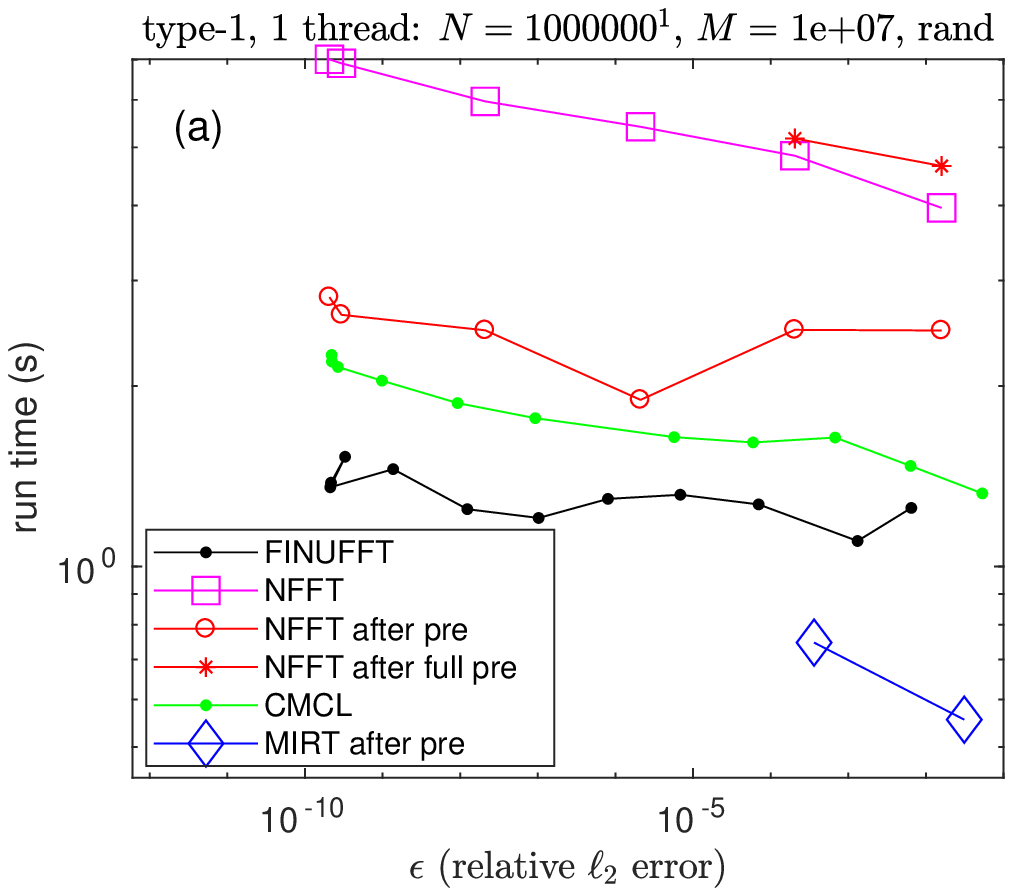}
\;\;
\ig{height=2.2in}{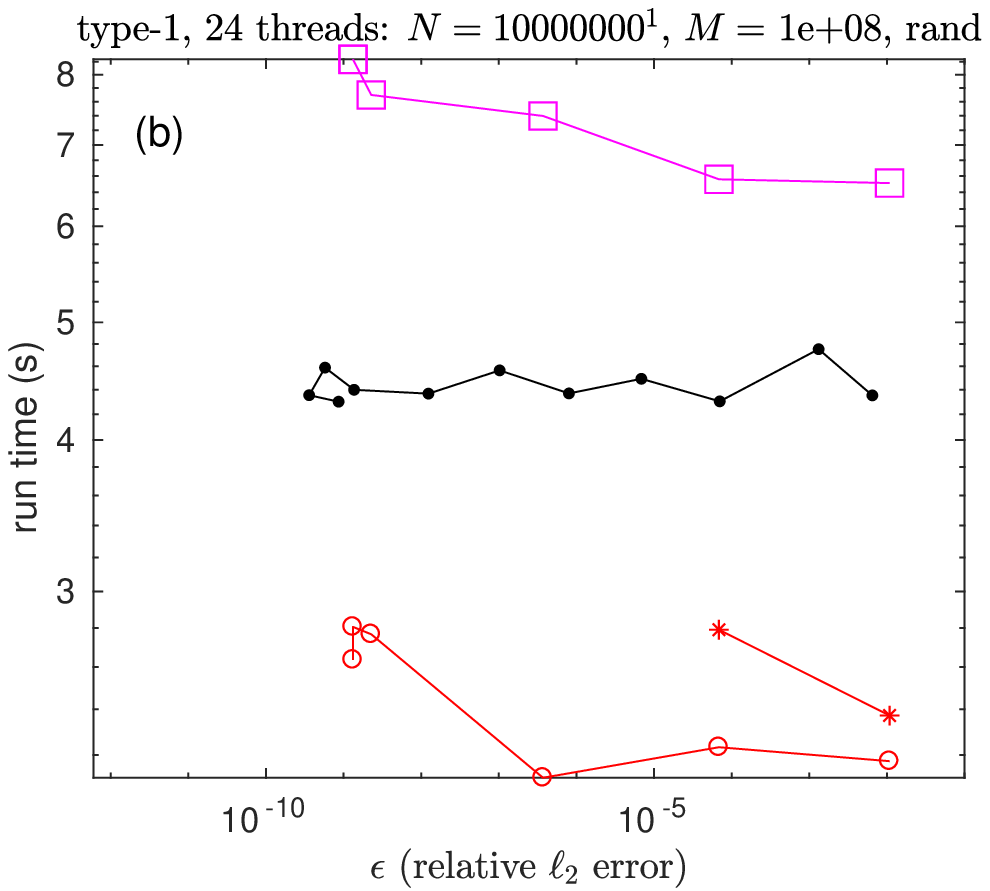}
\\
\vspace{0ex}\\
\mbox{}\;\;
\ig{height=2.2in}{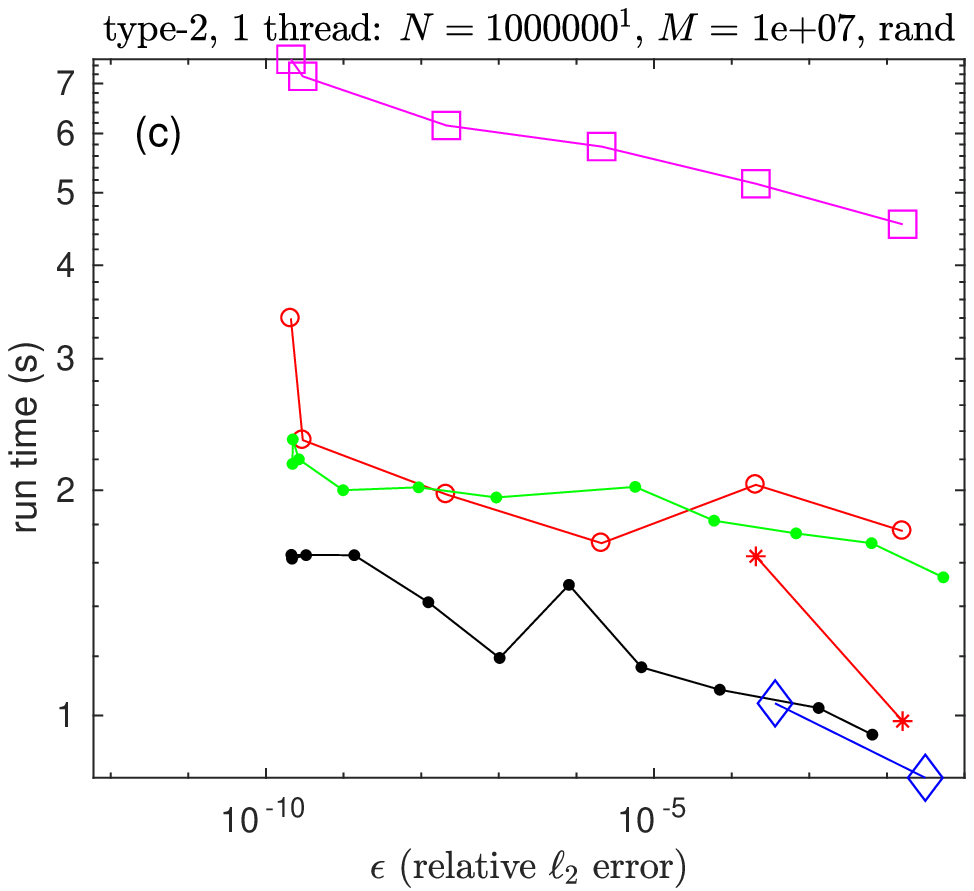}
\;\;
\ig{height=2.2in}{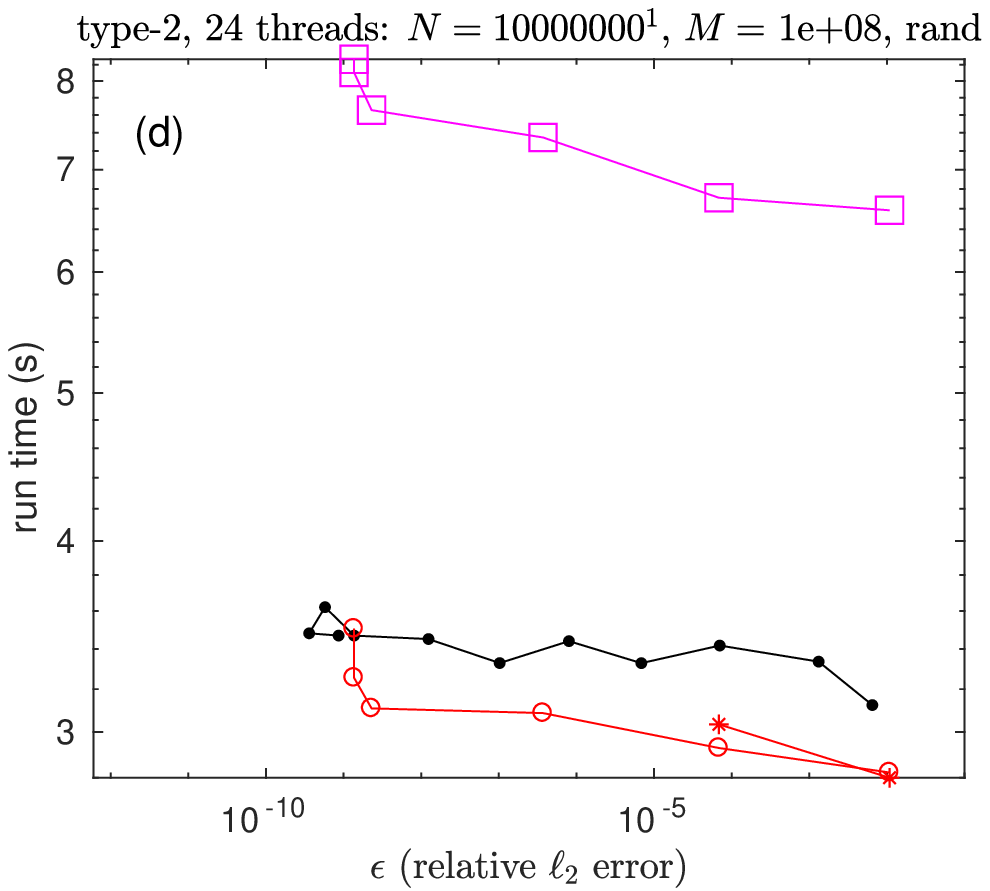}
\ca{1D comparisons.
  Execution time vs accuracy is shown for various NUFFT libraries,
  for random data in 1D.
  Precomputations (needed for codes labeled ``pre'') were not included.
  The left shows single-threaded codes, the right multi-threaded.
  The top pair are type 1, the bottom pair type 2.
  See section~\ref{s:bench}.
}{f:1d}
\efi

\subsection{Benchmarks against existing libraries}
\label{s:bench}

We now compare
FINUFFT against several popular open-source CPU-based NUFFT libraries
mentioned in section~\ref{s:prior}.
Their properties (including their periodic domain conventions)
are summarized in Table~\ref{t:codes}.
We study speed vs accuracy for types 1 and 2 for $d=1,2,3$,
covering most applications.
The machine, OS, and default compiler were as above.
In multi-threaded tests we set $p=24$ (two threads per core).
Each code provides a MATLAB interface, or is native MATLAB.
For reproducibility, we now list their test parameters and set-up
(also see {\tt https://github.com/ahbarnett/nufft-bench} ):
\bi
\item FINUFFT, version 1.0. Compiler flags are as in the previous section.
  We tested tolerances $\eps=10^{-2}, 10^{-3},\ldots,10^{-12}$.

\item CMCL NUFFT, version 1.3.3 \cite{cmcl}.
  This ships with MEX binaries dated 2014. It uses {\tt dfftpack} for
  FFTs.
  For fairness, we recompiled the relevant {\tt *.mexa64} binaries
  on the test machine using {\tt gfortran} with flags
  {\tt -fPIC -Ofast -funroll-loops -march=native}, and {\tt mex}.
   We tested tolerances $10^{-1}, 10^{-2},\ldots,10^{-11}$.
 
 \item NFFT, version 3.3.2 \cite{nfft}.
   A compiler error resulted with GCC 7.x, so we used GCC 6.4.0.
   We used the default (backward KB) kernel.
   We used the ``guru'' interface with FFT grid size set to the smallest
   power of two at least $2N_i$, where $N_i$ is the number of modes in dimension
   $i$, following their examples.
   Since they increased speed, we set the flags
   \verb+PRE_PHI_HUT+,
   \verb+FFT_OUT_OF_PLACE+,
   \verb+NFFT_OMP_BLOCKWISE_ADJOINT+,
   and \verb+NFFT_SORT_NODES+ (the latter is not part of the standard
   MATLAB interface).
   We tested three variants:
  \bi
\item no kernel precomputation; kernel is evaluated on the fly
  (labeled ``NFFT'').
\item ``pre'': option \verb+PRE_PSI+ which precomputes
  $wdM$ kernel values (with tensor products done on the fly).
\item ``full pre'': option \verb+PRE_FULL_PSI+ which precomputes
  and then looks up all $w^dM$ kernel values.
  \ei
  We tested kernel parameters
  $m=1,2,\ldots,6$ (kernel width is $w=2m+1$),
  apart from ``full pre'' where RAM constraints limited us to $m=1,2$.

\item MIRT (Michigan Image Reconstruction Toolbox), no version number; however the latest changes to the
  {\tt nufft} directory were on Dec 13, 2016 \cite{MIRT}.
  This native MATLAB code precomputes a sparse matrix with all kernel values.
  Its matrix-vector multiplication appears to be single-threaded,
  thus we place this library in the single-threaded category.  
  We use oversampling $\rat=2.0$ and the default kernel {\tt minmax:kb},
  which appears to be close to KB \eqref{KB}.
  RAM constraints limited us to
  testing width parameters $J=2,4$ (equivalent to our $w$).
  
\item BART (Berkeley Advanced Reconstruction Toolbox), version 0.4.02
  \cite{BART}.
  This is a recent multithreaded
  C code for 3D only by Uecker, Lustig and Ong, 
  used in a recent comparison by Ou \cite{ou}.
  We compiled with {\tt -O3 -ffast-math -funroll-loops -march=native}.
  The MATLAB interface writes to and reads from temporary files; however
  for our problem sizes with the use of a local SSD drive this adds
  less than 10\% to the run-time.
  BART did not ship with periodic wrapping of the spreading kernel;
  however, upon request\footnote{M.~Uecker, private communication.}
  we received a patched code {\tt src/noncart/grid.c}.
  It has fixed accuracy ($w=3$ is fixed).
  We empirically find that a prefactor $\sqrt{N_1N_2N_3}/1.00211$ gives around 5-digit accuracy (without the strange factor $1.00211$ it gives only 3 digits).

\ei
The above notes also illustrate some of the challenges in setting up fair
comparisons.

\bfi
\ig{height=2.2in}{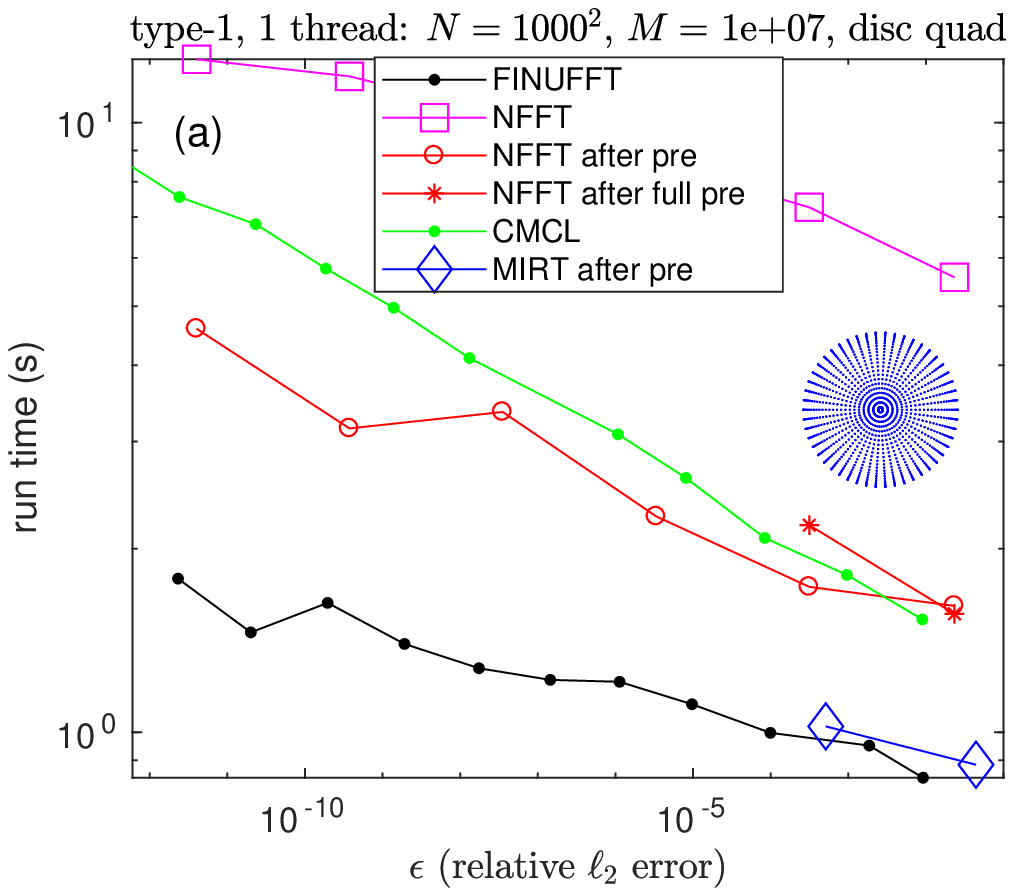}
\;\;
\ig{height=2.2in}{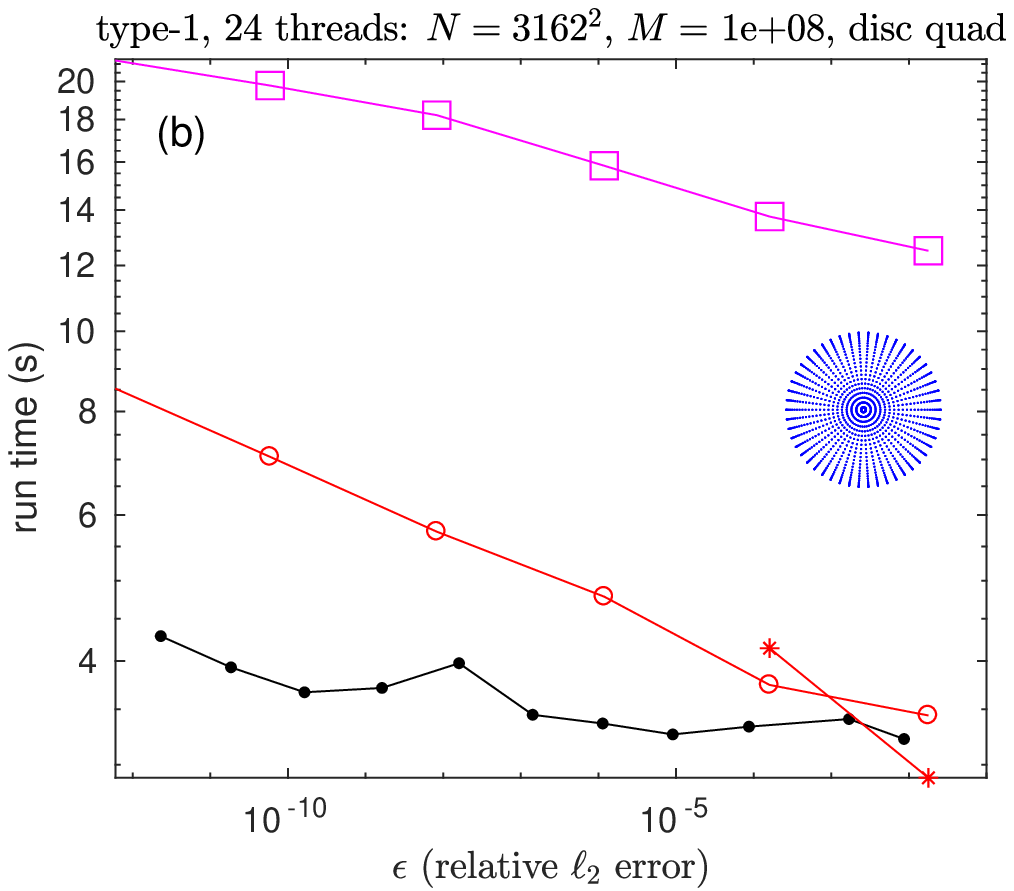}
\\
\vspace{0ex}\\
\ig{height=2.2in}{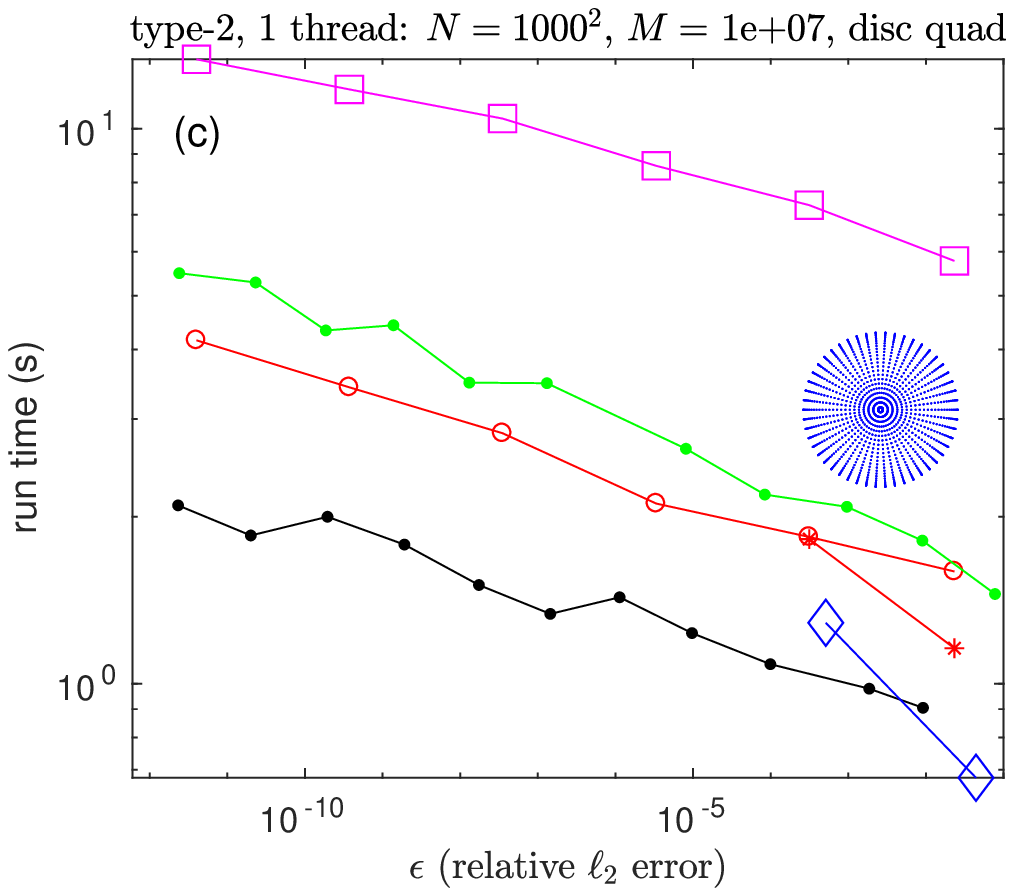}
\;\;
\ig{height=2.2in}{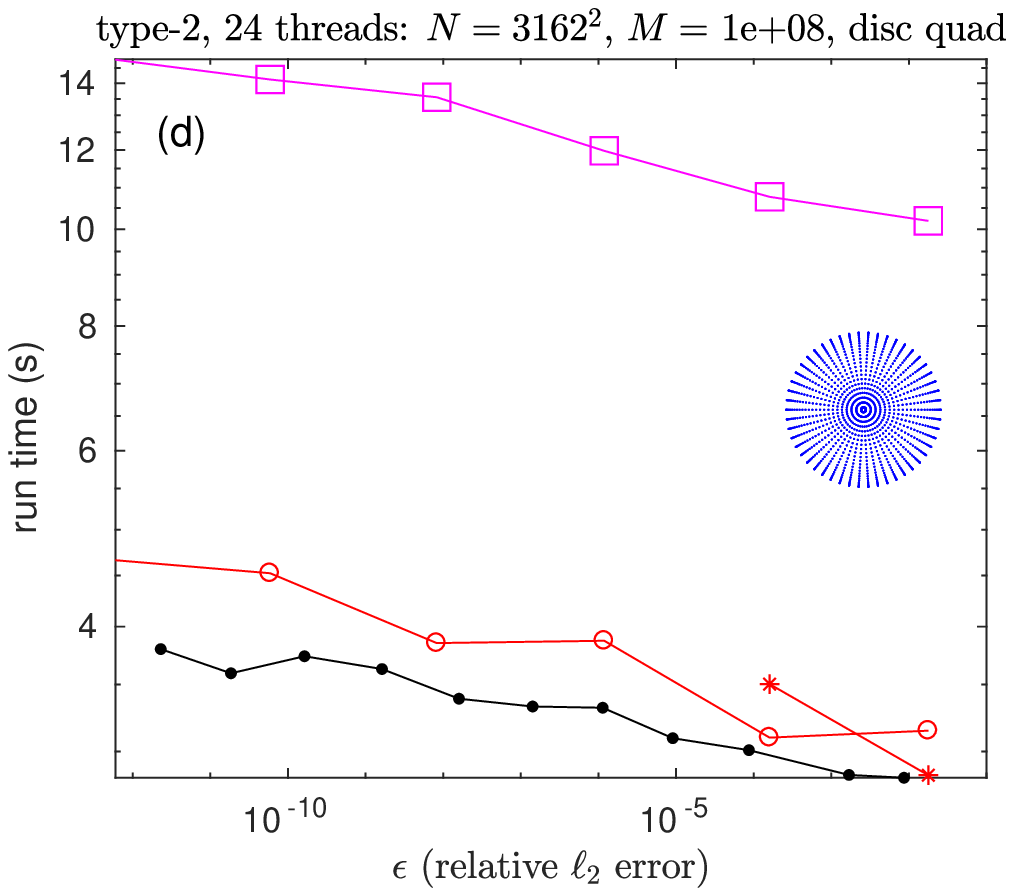}
\ca{2D comparisons.
  Execution time vs accuracy is shown for the tested libraries, for 2D polar
  ``disc quad'' nodes (see section~\ref{s:perf})
  illustrated in the insets at smaller $M$.
  Precomputations (needed for codes labeled ``pre'') were not included.
  The left shows single-threaded codes, the right multi-threaded.
  The top pair are type 1, the bottom pair type 2.
  See section~\ref{s:bench}.
}{f:2d}
\efi

We now discuss the results (Figs.~\ref{f:1d}--\ref{f:3dmulti}).
In each case we chose $M\approx 10 N$, for reasons discussed above.
To be favorable to codes that require precomputation,
precomputation times were not counted (hence the label
``after pre'' in the figures).
As in section~\ref{s:params}, $\epsilon$ denotes
relative $l_2$ error, measured against a ground truth of
FINUFFT with tolerance $\eps=10^{-14}$.

{\bf 1D comparisons.}\;
Fig.~\ref{f:1d} compares single-thread codes (left plots),
and then, for a larger task, multi-threaded codes (right plots).
For single-threaded, FINUFFT outperforms all libraries except
MIRT, which exploits MATLAB's apparently efficient sparse
matrix-vector multiplication.
However, what is not shown is that the precomputation time for
MIRT is around {\em 100 times longer} than the transform time.
For multithreaded type 1, FINUFFT is around $1.5$--$2\times$ faster
than NFFT without precomputation, but around $2\times$ slower than
``NFFT pre''.
For type 2, FINUFFT and ``NFFT pre'' are similar,
but of course this does not count the precomputation time
(and higher RAM overhead) of ``NFFT pre''.
As per Remark~\ref{r:rounding}, all libraries bottom out at around 9-10
digits due to rounding error.

{\bf 2D comparisons.}\;
Fig.~\ref{f:2d} shows similar comparisons in 2D, for 
a point distribution concentrating at the origin.
Compared to other codes not needing precomputation,
FINUFFT is 2--$5\times$ faster than CMCL (when single-threaded),
and 4--$8\times$ faster than NFFT.
When NFFT is allowed precomputation, its type 2 multi-threaded
speed is similar to FINUFFT, but for type 1 FINUFFT is $2\times$ faster
at high accuracy.

\bfi
\ig{height=2.2in}{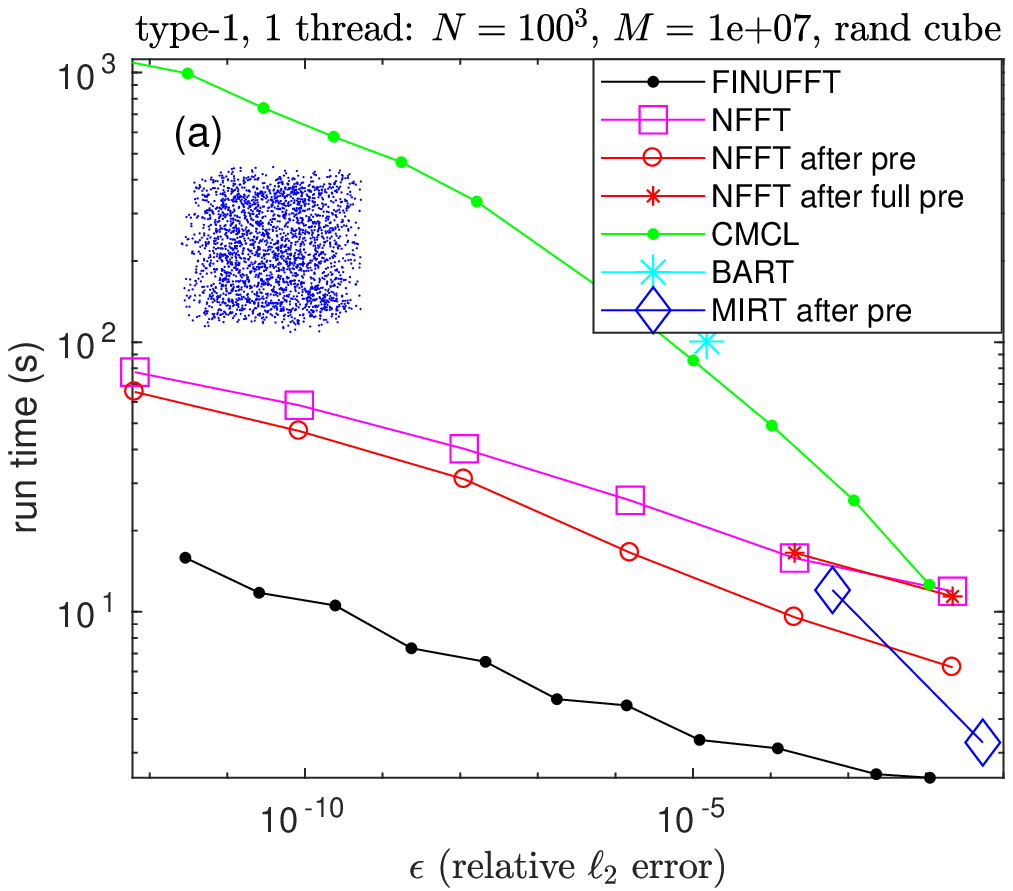}
\;\;
\ig{height=2.2in}{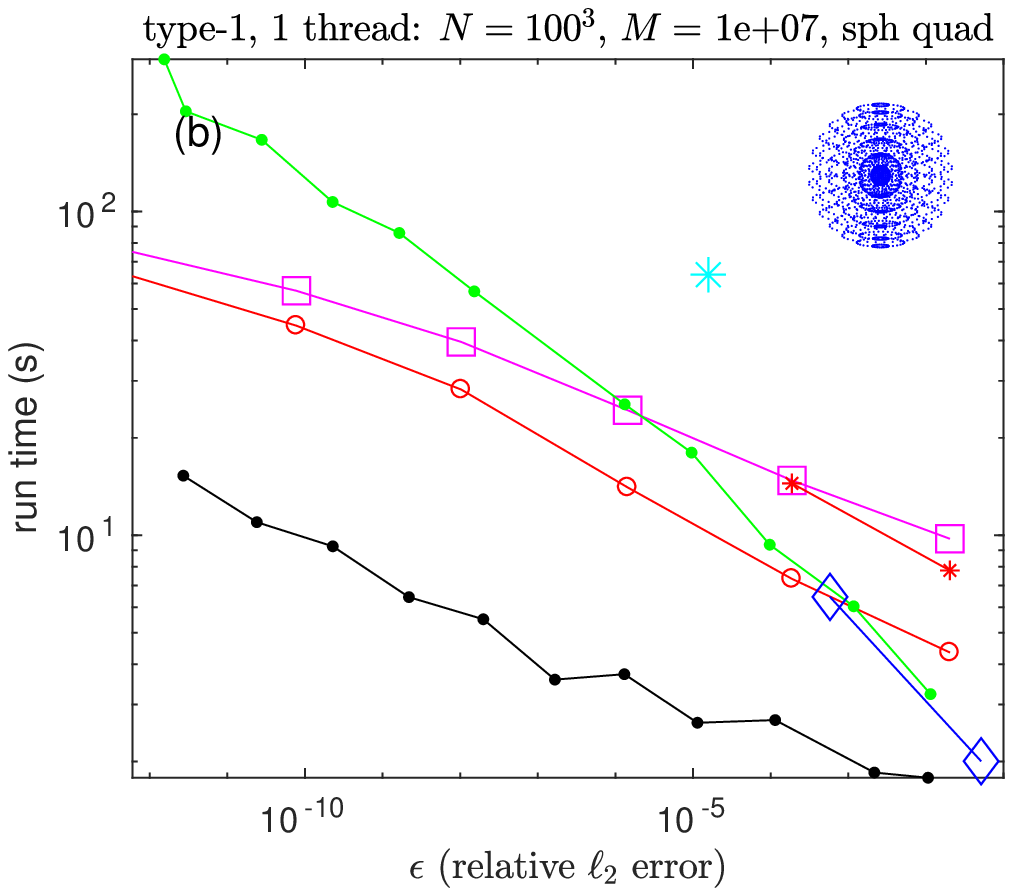}
\\
\vspace{0ex}\\
\ig{height=2.2in}{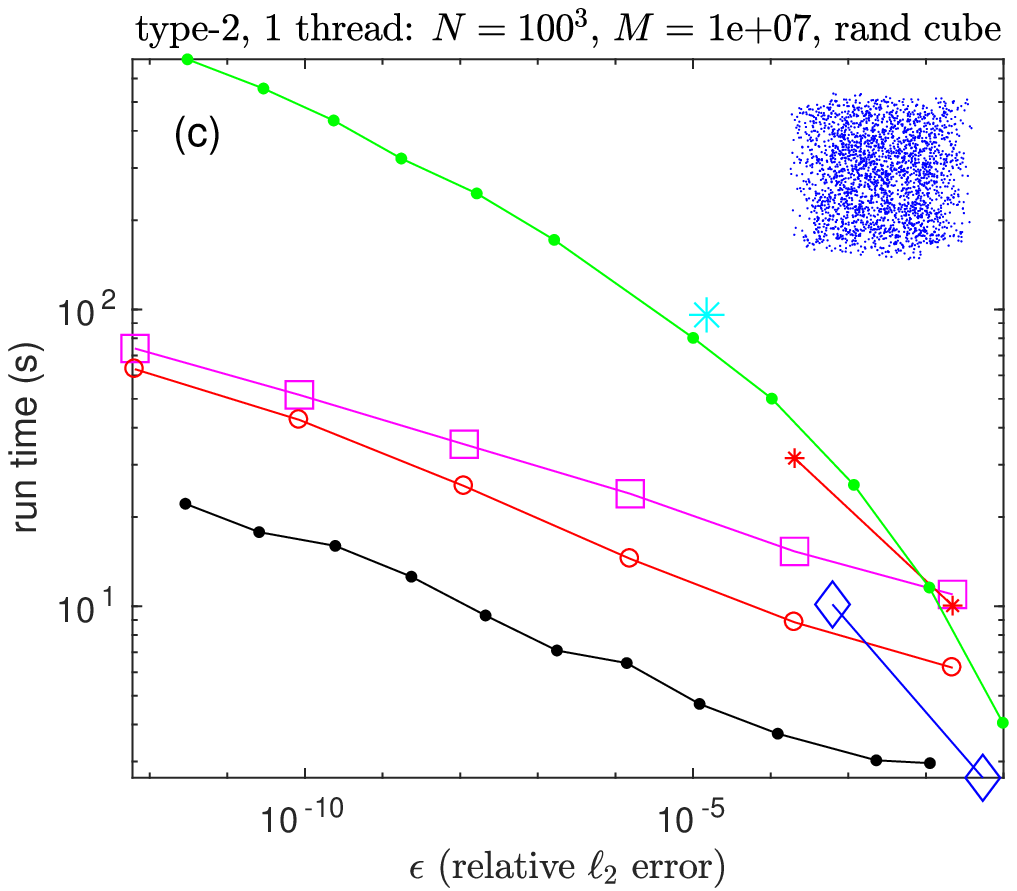}
\;\;
\ig{height=2.2in}{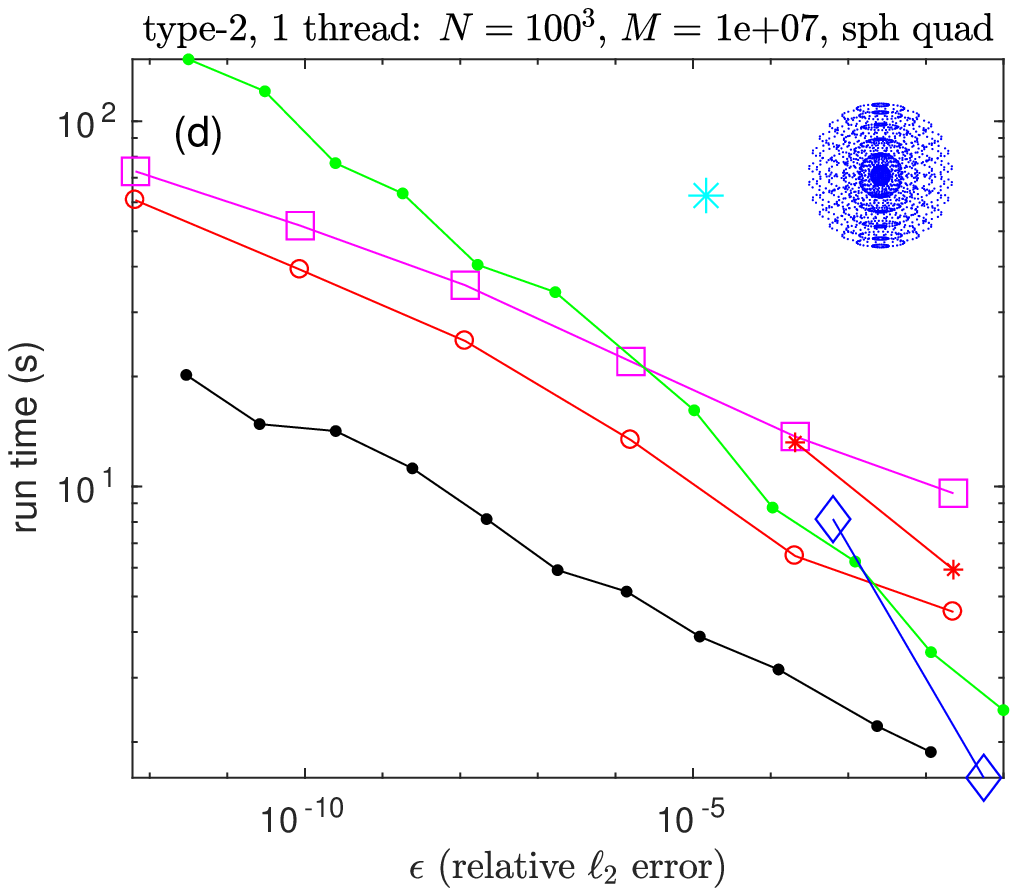}
\ca{3D single-threaded comparisons.
  Execution time vs accuracy is shown for the tested libraries,
  for uniform random (left pair) and spherical quadrature (right pair) nodes
  (see section~\ref{s:perf}).
  Node patterns are shown in the insets at a smaller $M$.
  Precomputations (needed for codes labeled ``pre'') were not included.
  The top pair are type 1, the bottom pair type 2.
  See section~\ref{s:bench}.
}{f:3dsingle}
\efi

\bfi
\ig{height=2.2in}{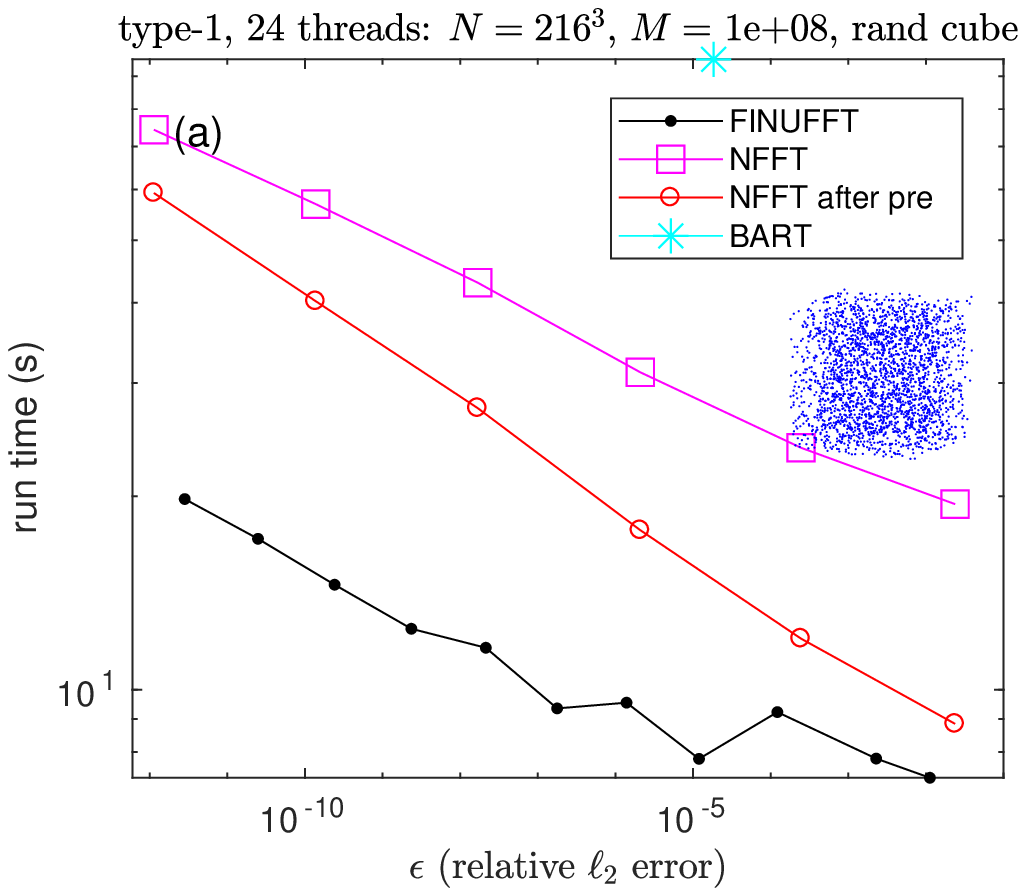}
\;\;
\ig{height=2.2in}{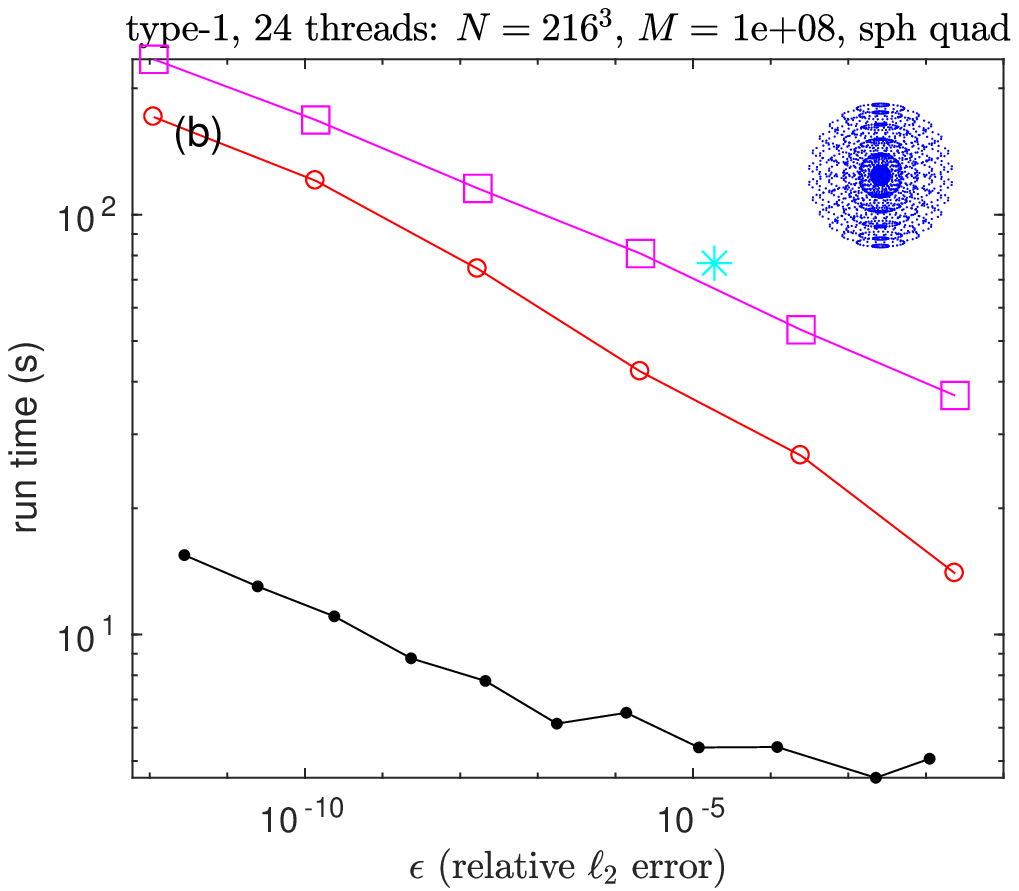}
\\
\vspace{0ex}\\
\ig{height=2.2in}{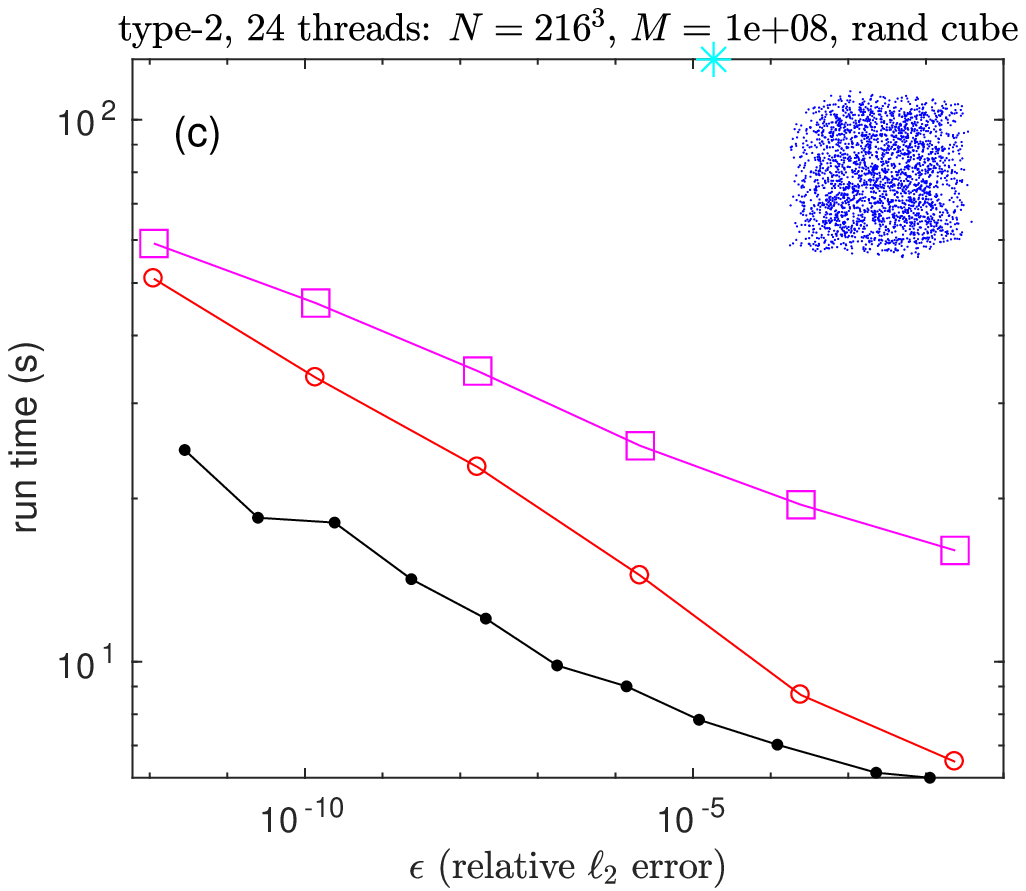}
\;\;
\ig{height=2.2in}{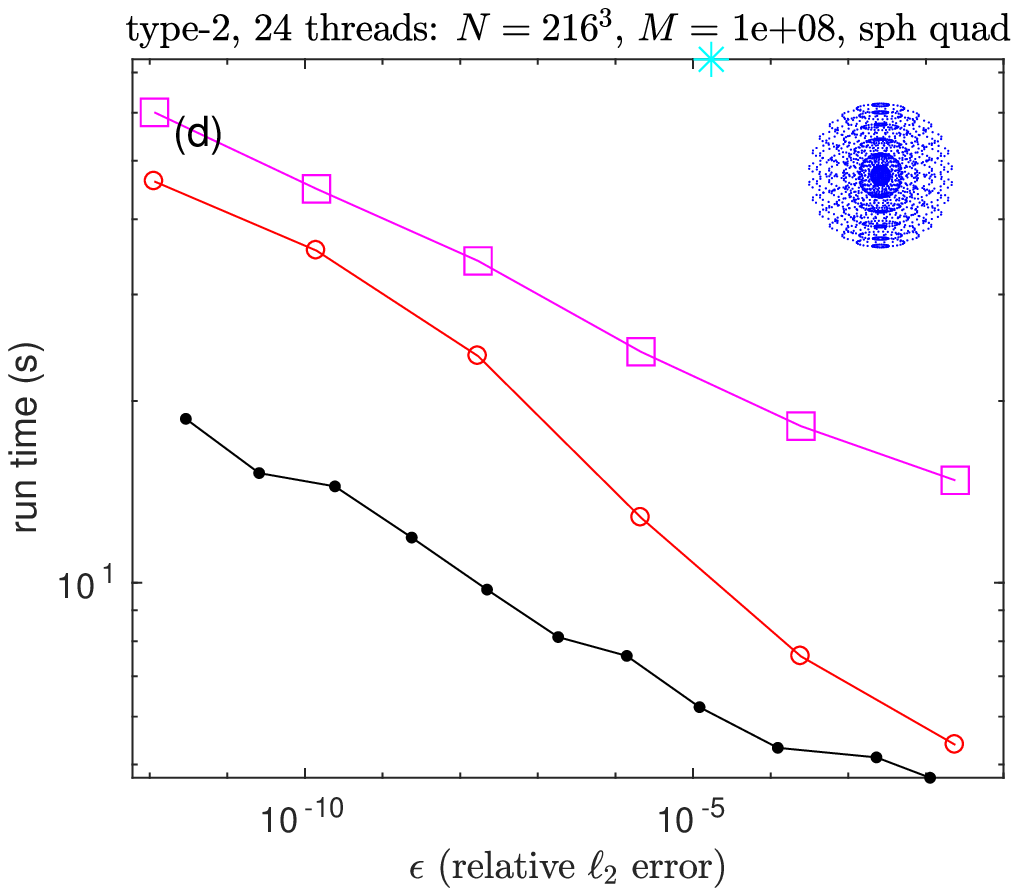}
\ca{3D multithreaded comparisons. For explanation see caption of
  Fig.~\ref{f:3dsingle}.
}{f:3dmulti}
\efi

\begin{table} 
\centering
\begin{tabular}{lllll}
    Code and parameters & $\epsilon$ (rel. $\ell_2$ error) & $t_\tbox{plan+pre}$ & $t_\tbox{run}$ & RAM overhead\\
    \hline
 FINUFFT (tol. $10^{-6}$) & 1.4e-06  &           N/A    &        14.6 s    &        8.8 GB\\
 NFFT ($m=3$) & 4.7e-06   &         10.4 s    &        238 s    &       20.9 GB\\
 NFFT ($m=3$) \verb+PRE_PSI+ & 4.7e-06   &         67.3 s    &        125 s    &       67.1 GB\\
 \hline
\end{tabular}
  \vspace{1ex}
  \ca{Performance of FINUFFT and NFFT for a large 3D Type 1 transform
    with roughly 6-digit accuracy, using 24 threads.
    $N=256^3$ modes are requested with $M=3\times10^8$ ``sph quad''
    \NU\ points.
    The spreading time dominates over the FFT.
    For NFFT, $t_\tbox{plan+pre}$ counts both
    planning and kernel precomputation.
    RAM is measured using {\tt top}, relative to the baseline
    (around 12 GB) needed to store the input data in MATLAB.
    See section~\ref{s:perf} for machine and NFFT parameters.
  }{t:big}
\end{table} 

{\bf 3D comparisons.}\;
Fig.~\ref{f:3dsingle} compares single-threaded codes (now including BART).
The left pair of plots shows
random points: FINUFFT is at least $2\times$ faster than any other code,
apart from MIRT at 1-digit accuracy.
CMCL is a factor 4--$50\times$ slower than single-threaded FINUFFT,
we believe in part due to its lack of sorting \NU\ points.
(The evidence is that for the right pair,
where points have an ordered access pattern,
CMCL is only 2--$10\times$ slower).
NFFT without precomputation is 3--$5\times$ slower than FINUFFT;
precomputation brings this down to 2--$4\times$.
As for $d=1,2$, we observe that ``NFFT full pre'' is no faster than
``NFFT pre'', despite its longer precomputation and larger RAM overhead.

Fig.~\ref{f:3dmulti} shows larger multi-threaded comparisons against NFFT;
now we cannot include ``NFFT full pre'' due to its large RAM usage.
At low accuracies with random points, FINUFFT and ``NFFT pre'' have similar
speeds.
However, for type 1 ``sph quad'' (panel (b)),
for $\epsilon<10^{-6}$, FINUFFT is 8--$10\times$ faster than NFFT
even with precomputation.
FINUFFT is at least $10\times$ faster than BART in all cases.

In Table~\ref{t:big} we
compare FINUFFT and NFFT in terms of both speed and memory overhead,
for the same large, medium-accuracy, multi-threaded task.
We emphasize that, since $N_i=256$, $i=1,2,3$, is a power of two,
the fine grids chosen by the two libraries are an identical $n_i=512$.
This means that the memory use, and the FFTW calls, are identical.
Furthermore, the kernel widths are both $w=7$ so the numbers of fine grid
points written to are identical.
If precomputations are excluded, FINUFFT is
$16\times$ faster than NFFT, and $8.6\times$ faster than ``NFFT pre''.
For a single use (\ie including initialization and precomputation),
these ratios become $17\times$ and $13\times$, respectively.

\begin{rmk}
  We believe that the following explains the large type 1 performance gain of FINUFFT over NFFT
  for the 3D clustered ``sph quad'' nodes,
  shown by Fig.~\ref{f:3dmulti}(b) and Table~\ref{t:big}.
  NFFT assigns threads to equal slices
  of the fine grid (in the $x$-direction),
  whereas FINUFFT uses subproblems which load balance regardless of the
  clustering of nodes.
  We find that only a couple of threads are active for the entire
  run time of NFFT; most complete their jobs quickly,
  giving low parallel efficiency.
  \label{r:loadbal}
\end{rmk}

Finally, Table~\ref{t:big} shows that, if NFFT precomputation is used,
its RAM overhead is around $8\times$ that of FINUFFT.
The ``NFFT pre'' RAM overhead of around 28 {\tt double}s per point
is consistent with the expected $wd=21$ stored kernel values
per point.

\section{Conclusion}
\label{s:conc}

We have presented an open-source parallel CPU-based general purpose
type 1, 2 and 3
NUFFT library (FINUFFT)
for dimensions 1, 2, and 3
that is competitive with existing CPU-based libraries.
Using a new spreading kernel \eqref{ES}, all kernel
evaluations are efficiently done on the fly:
this avoids any precomputation phase, keeps the RAM overhead small,
and allows for a simple user-friendly interface from multiple languages.
Efficient parallelization balances the work of threads, adapting to
any \NU\ point distribution.
For all three types, we introduce numerical quadrature to evaluate a
kernel \FT\ for which there is no known analytic formula.
Rigorous estimates show almost exponential convergence of the
kernel aliasing error, with rate arbitrarily close to that of the best known.
We explained the gap between such estimates and empirical relative errors.
We benchmarked several NUFFT libraries in detail.
We showed that for certain 3D problems with clustered distributions
FINUFFT is an order of magnitude faster than
the other libraries, even when they are allowed precomputation.
In the latter case, FINUFFT has an order of magnitude less RAM overhead.

There are several directions for future work,
starting with benchmarking the less-common type 3 case.
An efficient interface for the case of repeated small problems
(Remark~\ref{r:small}) should be completed.
There is also a need for a carefully-benchmarked
general-purpose GPU NUFFT code, following
Kunis--Kunis \cite{cunfft}, Ou \cite{ou}, and others.
Implementation of both of the above is in progress.


\begin{rmk}
  In some particle-mesh Ewald applications \cite{Lindbo11}
  one needs spatial derivatives of the spreading kernel.%
  \footnote{A.-K. Tornberg, personal communication.}
  However, \eqref{ES} has unbounded derivatives (with inverse square-root
  singularity) at the endpoints.
  Instead one may prefer the exponentially-close
  variant $\phi(z) = 2e^{-\freq} \cosh \freq \sqrt{1-z^2}$
  since it is smooth up to the endpoints.
  This kernel requires one extra reciprocal, or approximation
  as in section~\ref{s:poly}.
\end{rmk}

\section*{Acknowledgments}

We are grateful for discussions with
Joakim And\'en,
Charlie Epstein,
Zydrunas Gimbutas,
Leslie Greengard,
Hannah Lawrence,
Andras Pataki,
Daniel Potts,
Vladimir Rokhlin,
Yu-hsuan Shih,
Marina Spivak,
David Stein,
and
Anna-Karin Tornberg.
The Flatiron Institute is a division of the Simons Foundation.

\bibliographystyle{abbrv}
\bibliography{finufft_re}

\begin{thebibliography}{10}

\bibitem{finufft.jl}
L.~af~Klinteberg.
\newblock {J}ulia interface to {FINUFFT}, 2018.
\newblock {\tt https://github.com/ludvigak/FINUFFT.jl}.

\bibitem{anderson96}
C.~Anderson and M.~D. Dahleh.
\newblock Rapid computation of the discrete {F}ourier transform.
\newblock {\em SIAM J. Sci. Comput.}, 17(4):913--919, 1996.

\bibitem{andersson12}
F.~Andersson, R.~Moses, and F.~Natterer.
\newblock Fast {F}ourier methods for synthetic aperture radar imaging.
\newblock {\em IEEE Trans. Aerospace Elec. Sys.}, 48(1):215--229, 2012.

\bibitem{nufftanal}
A.~H. Barnett.
\newblock Asymptotic aliasing error of the kernel exp$(\beta \sqrt{1-x^2})$
  when used for non-uniform fast {F}ourier transforms, 2019.
\newblock in preparation.

\bibitem{finufftgit}
A.~H. Barnett, J.~Magland, and L.~af~Klinterberg.
\newblock {F}latiron {I}nstitute nonuniform fast {F}ourier transform libraries
  ({FINUFFT}), 2018.
\newblock {\tt https://github.com/flatironinstitute/finufft}.

\bibitem{cryo}
A.~H. Barnett, M.~Spivak, A.~Pataki, and L.~Greengard.
\newblock Rapid solution of the cryo-{EM} reconstruction problem by frequency
  marching.
\newblock {\em SIAM J. Imaging Sci.}, 10(3):1170--1195, 2017.

\bibitem{beylkinnufft}
G.~Beylkin.
\newblock On the fast {F}ourier transform of functions with singularities.
\newblock {\em Appl. Comput. Harmonic Anal.}, 2:363--383, 1995.

\bibitem{boettcher}
A.~B\"ottcher and D.~Potts.
\newblock Probability against condition number and sampling of multivariate
  trigonometric random polynomials.
\newblock {\em Electron. Trans. Numer. Anal.}, 26:178--189, 2007.

\bibitem{ultranufft}
M.~M. Bronstein, A.~M. Bronstein, M.~Zibulevsky, and H.~Azhari.
\newblock Reconstruction in diffraction ultrasound tomography using nonuniform
  {FFT}.
\newblock {\em IEEE Trans.\ Medical Imaging}, 21(11):1395--1401, 2002.

\bibitem{usingopenmp}
B.~Chapman, G.~Jost, and R.~van~der Pas.
\newblock {\em Using {OpenMP}. Portable Shared Memory Parallel Programming}.
\newblock MIT Press, 2008.

\bibitem{davisheller}
M.~J. Davis and E.~J. Heller.
\newblock Semiclassical {G}aussian basis set method for molecular vibrational
  wave functions.
\newblock {\em J. Chem. Phys.}, 71(8):3383--3395, 1979.

\bibitem{dutt}
A.~Dutt and V.~Rokhlin.
\newblock Fast {F}ourier transforms for nonequispaced data.
\newblock {\em SIAM J. Sci. Comput.}, 14:1369--1393, 1993.

\bibitem{duttcotfmm}
A.~Dutt and V.~Rokhlin.
\newblock Fast {F}ourier transforms for nonequispaced data, {II}.
\newblock {\em Appl. Comput. Harmonic Anal.}, 2:85--100, 1995.

\bibitem{elbel}
B.~Elbel and G.~Steidl.
\newblock Fast {F}ourier transform for nonequispaced data.
\newblock In {\em Approximation Theory {IX}}, pages 39--46. Vanderbilt
  University Press, Nashville, 1998.

\bibitem{MIRT}
J.~Fessler.
\newblock Michigan image reconstruction toolbox, 2016.
\newblock {\tt https://web.eecs.umich.edu/$\sim$fessler/irt/fessler.tgz}.

\bibitem{fessler}
J.~Fessler and B.~Sutton.
\newblock Nonuniform fast fourier transforms using min-max interpolation.
\newblock {\em IEEE Trans.\ Signal Proc.}, 51(2):560--574, 2003.

\bibitem{fourmontthesis}
K.~Fourmont.
\newblock {\em Schnelle {F}ourier-{T}ransformation bei nicht\"aquidistanten
  {G}ittern und tomographische {A}nwendungen}.
\newblock PhD thesis, Univ. M\"unster, 1999.

\bibitem{fourmont}
K.~Fourmont.
\newblock Non-equispaced fast {F}ourier transforms with applications to
  tomography.
\newblock {\em J. Fourier Anal. Appl.}, 9(5):431--450, 2003.

\bibitem{fftw}
M.~Frigo and S.~G. Johnson.
\newblock {FFTW}.
\newblock {\tt http://www.fftw.org/}.

\bibitem{gimbutasgrid}
Z.~Gimbutas and S.~Veerapaneni.
\newblock A fast algorithm for spherical grid rotations and its application to
  singular quadrature.
\newblock {\em SIAM J. Sci. Comput.}, 5(6):A2738--A2751, 2013.

\bibitem{kaiserinterview}
A.~Goldstein and J.~Abbate.
\newblock Oral history: {J}ames {K}aiser.
\newblock {\tt http://ethw.org/Oral-History:James\_Kaiser}, 1997.
\newblock online; accessed 2017-04-15.

\bibitem{GS8}
I.~Gradshteyn and I.~Ryzhik.
\newblock {\em Table of Integrals, Series and Products}.
\newblock New York: Academic, 8th edition, 2015.

\bibitem{fastsinc}
L.~Greegard, J.-Y. Lee, and S.~Inati.
\newblock The fast sinc transform and image reconstruction from nonuniform
  samples in $k$-space.
\newblock {\em Comm. Appl. Math. and Comp. Sci.}, 1(1):121--131, 2006.

\bibitem{cmcl}
L.~Greengard and J.-Y. Lee.
\newblock {NUFFT} libraries in {F}ortran.
\newblock {\tt http://www.cims.nyu.edu/cmcl/nufft/nufft.html}.
\newblock online; accessed 2017-04-10.

\bibitem{nufft}
L.~Greengard and J.-Y. Lee.
\newblock Accelerating the nonuniform fast {F}ourier transform.
\newblock {\em SIAM Review}, 46(3):443--454, 2004.

\bibitem{jackson91}
J.~I. Jackson, C.~H. Meyer, D.~G. Nishimura, and A.~Macovski.
\newblock Selection of a convolution function for {F}ourier inversion using
  gridding.
\newblock {\em IEEE Trans.\ Medical Imaging}, 10(3):473--478, 1991.

\bibitem{l2jacob}
M.~Jacob.
\newblock Optimized least-square nonuniform fast {F}ourier transform.
\newblock {\em IEEE Trans. Signal Process.}, 57(6):2165--2177, 2009.

\bibitem{kaiser}
J.~Kaiser.
\newblock Digital filters.
\newblock In J.~Kaiser and F.~Kuo, editors, {\em System analysis by digital
  computer}, chapter~7, pages 218--285. Wiley, 1966.

\bibitem{nfft}
J.~Keiner, S.~Kunis, and D.~Potts.
\newblock {NFFT}.
\newblock {\tt http://www-user.tu-chemnitz.de/$\sim$potts/nfft}, 2002--2016.
\newblock online; accessed 2017-04-10.

\bibitem{usingnfft}
J.~Keiner, S.~Kunis, and D.~Potts.
\newblock Using {NFFT} 3 --- a software library for various nonequispaced fast
  {F}ourier transforms.
\newblock {\em ACM Trans. Math. Software}, 36(4), 2009.

\bibitem{kircheis}
M.~Kircheis and D.~Potts.
\newblock Direct inversion of the nonequispaced fast fourier transform, 2019.
\newblock {\tt arxiv:1811.05335}, submitted to {\em Linear Algebra Appl.}

\bibitem{gpunufft}
F.~Knoll, A.~Schwarzl, C.~Diwoki, and D.~K. Sodickson.
\newblock {gpuNUFFT}---an open-source {GPU} library for {3D} gridding with
  direct {MATLAB} interface.
\newblock {\em Proc. Intl. Soc. Mag. Reson. Med.}, 22, 2014.
\newblock {\tt https://github.com/andyschwarzl/gpuNUFFT}.

\bibitem{kunisthesis}
S.~Kunis.
\newblock {\em Nonequispaced {FFT}: Generalisation and Inversion}.
\newblock PhD thesis, Universit\"at zu L\"ubeck, 2006.

\bibitem{cunfft}
S.~Kunis and S.~Kunis.
\newblock The nonequispaced {FFT} on graphics processing units.
\newblock {\em Proc. Appl. Math. Mech.}, 12:7--10, 2012.

\bibitem{nufft3}
J.-Y. Lee and L.~Greengard.
\newblock The type 3 nonuniform {FFT} and its applications.
\newblock {\em J. Comput. Phys.}, 206:1--5, 2005.

\bibitem{pynufft}
J.-M. Lin.
\newblock Python non-uniform fast {F}ourier transform ({PyNUFFT}):
  multi-dimensional non-{C}artesian image reconstruction package for
  heterogeneous platforms and applications to {MRI}, 2017.
\newblock {\tt arxiv:1710.03197v1}.

\bibitem{Lindbo11}
D.~Lindbo and A.-K. Tornberg.
\newblock Spectral accuracy in fast {E}wald-based methods for particle
  simulations.
\newblock {\em J. Comput. Phys.}, 230:8744--8761, 2011.

\bibitem{emnufft}
Q.~H. Liu and N.~Nguyen.
\newblock An accurate algorithm for nonuniform fast fourier transforms
  ({NUFFT}'s).
\newblock {\em IEEE Microwave and Guided Wave Letters}, 8(1):18--20, 1998.

\bibitem{meisel78}
D.~D. Meisel.
\newblock {F}ourier transforms of data samples at unequal observational
  intervals.
\newblock {\em Astronom. J.}, 83(5):538--545, 1978.

\bibitem{CEPHES}
S.~L. Moshier.
\newblock {CEPHES} mathematical function library, 1984--1992.
\newblock {\tt http://www.netlib.org/cephes}.

\bibitem{nestler}
F.~Nestler.
\newblock Automated parameter tuning based on {RMS} errors for nonequispaced
  {FFT}s.
\newblock {\em Adv. Comput. Math.}, 42:889--919, 2016.

\bibitem{nestlerPME}
F.~Nestler, M.~Pippig, and D.~Potts.
\newblock Fast {E}wald summation based on {NFFT} with mixed periodicity.
\newblock {\em J. Comput. Phys.}, 285:280--315, 2015.

\bibitem{dlmf}
F.~W.~J. Olver, D.~W. Lozier, R.~F. Boisvert, and C.~W. Clark, editors.
\newblock {\em {NIST} Handbook of Mathematical Functions}.
\newblock Cambridge University Press, 2010.
\newblock {\tt http://dlmf.nist.gov}.

\bibitem{oppenheim71}
A.~Oppenheim, D.~Johnson, and K.~Steiglitz.
\newblock Computation of spectra with unequal resolution using the fast
  {F}ourier transform.
\newblock {\em Proc.\ IEEE}, 59(2):299--301, 1971.

\bibitem{osipov}
A.~Osipov, V.~Rokhlin, and H.~Xiao.
\newblock {\em Prolate Spheroidal Wave Functions of Order Zero: Mathematical
  Tools for Bandlimited Approximation}, volume 187 of {\em Applied Mathematical
  Sciences}.
\newblock Springer, US, 2013.

\bibitem{ou}
T.~Ou.
\newblock {gNUFFTW}: Auto-tuning for high-performance {GPU}-accelerated
  non-uniform fast {F}ourier transforms.
\newblock {\tt
  http://www2.eecs.berkeley.edu/Pubs/TechRpts/2017/EECS-2017-90.html}, 2017.
\newblock Technical Report \#UCB/EECS-2017-90, UC Berkeley.

\bibitem{pottshabil}
D.~Potts.
\newblock Schnelle {F}ourier-{T}ransformationen f{\"u}r nicht{\"a}quidistante
  {D}aten und {A}nwendungen, 2003.
\newblock Habilitationssschift, L{\"u}beck.

\bibitem{nfftchap}
D.~Potts, G.~Steidl, and M.~Tasche.
\newblock Fast {F}ourier transforms for nonequispaced data: a tutorial.
\newblock In {\em Modern Sampling Theory: Mathematics and Applications}, pages
  247--270. Birkh\"auser, Boston, 2001.

\bibitem{rybicki}
W.~H. Press and G.~B. Rybicki.
\newblock Fast algorithm for spectral analysis of unevenly sampled data.
\newblock {\em Astrophys. J.}, 338:277--280, 1989.

\bibitem{numrec}
W.~H. Press, S.~A. Teukolsky, W.~T. Vetterling, and B.~P. Flannery.
\newblock {\em Numerical recipes in {C}}.
\newblock Cambridge University Press, Cambridge, 2nd edition, 2002.

\bibitem{prudnikov1}
A.~Prudnikov, Y.~A. Brychkov, and O.~I. Marichev.
\newblock {\em Integrals and series, Volume 1. Elementary functions}.
\newblock Gordon and Breach, 1986.

\bibitem{prudnikov2}
A.~Prudnikov, Y.~A. Brychkov, and O.~I. Marichev.
\newblock {\em Integrals and series, Volume 2. Special functions}.
\newblock Gordon and Breach, 1986.

\bibitem{townsendnufft}
D.~Ruis-Antol{\'i}n and A.~Townsend.
\newblock A nonuniform fast {Fourier} transform based on low rank
  approximation, 2018.

\bibitem{slepian65}
D.~Slepian.
\newblock Some asymptotic expansions for prolate spheroidal wave functions.
\newblock {\em Journal of Mathematics and Physics}, 44(1-4):99--140, 1965.

\bibitem{slepianrev}
D.~Slepian.
\newblock Some comments on {F}ourier analysis, uncertainty, and modeling.
\newblock {\em SIAM Review}, 25:379--393, 1983.

\bibitem{steidl98}
G.~Steidl.
\newblock A note on fast {F}ourier transforms for nonequispaced grids.
\newblock {\em Adv. Comput. Math.}, 9:337--352, 1998.

\bibitem{suttonfield}
B.~P. Sutton, D.~C. Noll, and J.~A. Fessler.
\newblock Fast, iterative image reconstruction for {MRI} in the presence of
  field inhomogeneities.
\newblock {\em IEEE Trans.\ Med.\ Imaging}, 22(2):178--188, 2003.

\bibitem{thompson74}
A.~R. Thompson and R.~N. Bracewell.
\newblock Interpolation and {F}ourier transformation of fringe visibilities.
\newblock {\em Astronom. J.}, 79:11--24, 1974.

\bibitem{BART}
M.~Uecker and M.~Lustig.
\newblock {BART} toolbox for computational magnetic resonance imaging, 2016.
\newblock DOI: 10.5281/zenodo.592960. Available at {\tt
  https://mrirecon.github.io/bart/}.

\bibitem{gelbrecon}
A.~Viswanathan, A.~Gelb, D.~Cochran, and R.~Renaut.
\newblock On reconstruction from non-uniform spectral data.
\newblock {\em J. Sci. Comput.}, 45(1):487--513, 2010.

\bibitem{volkmer}
T.~Volkmer.
\newblock {OpenMP} parallelization in the {NFFT} software library, 2012.
\newblock Preprint 2012-07, Faculty of Mathematics, Technische Universit\"at
  Chemnitz.

\bibitem{warenufft}
A.~F. Ware.
\newblock Fast approximate {F}ourier transforms for irregularly spaced data.
\newblock {\em SIAM Review}, 40:838--856, 1998.

\bibitem{octnufft}
K.~Zhang and J.~U. Kang.
\newblock Graphics processing unit accelerated non-uniform fast {F}ourier
  transform for ultrahigh-speed, real-time {F}ourier-domain {OCT}.
\newblock {\em Opt. Express}, 18(22):23472--87, 2010.

\bibitem{steerablePCA}
Z.~Zhao, Y.~Shkolnisky, and A.~Singer.
\newblock Fast steerable principal component analysis.
\newblock {\em IEEE Trans.\ Comput.\ Imaging}, 2(1):1--12, 2016.

\end{thebibliography}
\end{document}